\RequirePackage{amsthm}

\documentclass[sn-mathphys-num]{sn-jnl}

\usepackage{adjustbox}
\usepackage[graphicx]{realboxes}
\usepackage{amsfonts}
\usepackage{geometry}
\usepackage{fancyhdr}
\usepackage{fancybox}
\usepackage{amsmath,amsthm,amssymb}
\usepackage{graphicx}
\usepackage[utf8]{inputenc}
\usepackage{enumerate}  
\usepackage{scrextend}
\usepackage{wasysym}
\usepackage{listings}
\usepackage{accents}
\usepackage{tikz}
\usepackage{float}
\usepackage{vwcol}
\usepackage{calc}
\usepackage{expl3}
\usepackage{parallel}
\usepackage{multirow}
\usepackage{tabularx}
\usepackage{xcolor,soul}
\usepackage{lipsum}
\usepackage{booktabs}
\usepackage{lmodern}
\usepackage{mathtools}
\usepackage{pdflscape}
\usepackage{autobreak}
\usepackage{empheq}
\usepackage{xstring}
\usepackage{epstopdf} 
\usepackage{scalerel}
\usepackage[percent]{overpic}
\usepackage{blindtext}
\usepackage{xfp}
\usepackage[title]{appendix}
\usepackage{mathrsfs}
\usepackage{textcomp}
\usepackage{manyfoot}
\usepackage{algorithm}
\usepackage{algorithmicx}
\usepackage{algpseudocode}
\usepackage{ulem}
\usepackage[font=small,labelfont=bf]{caption}
\usepackage{subcaption}
\usepackage{setspace}
\usepackage{CustomCommands-private}
\usepackage{orcidlink}

\usetikzlibrary{calc}
\usetikzlibrary{decorations.pathmorphing}
\usetikzlibrary{tikzmark,arrows,shapes,decorations.pathreplacing}
\tikzset{every picture/.style={remember picture}}

\DeclareFontFamily{OT1}{pzc}{}
\DeclareFontShape{OT1}{pzc}{m}{it}{<-> s * [1.10] pzcmi7t}{}
\DeclareMathAlphabet{\mathpzc}{OT1}{pzc}{m}{it}

\newcolumntype{P}[1]{>{\centering\arraybackslash}p{#1}}

\definecolor{TUMblue}{RGB}{0,101,189}
\definecolor{TUMgreen}{RGB}{162,173,0}
\definecolor{TUMdarkgreen}{RGB}{145,172,107}
\definecolor{TUMivory}{RGB}{218,215,203}
\definecolor{TUMorange}{RGB}{227,114,34}
\definecolor{TUMdarkred}{RGB}{202,033,063}
\definecolor{TUMkirschgruen}{RGB}{229,052,024}
\definecolor{LMUgreen}{RGB}{0,135,58} 
\hypersetup{colorlinks,linkcolor={TUMblue},citecolor={TUMdarkred},urlcolor=TUMdarkgreen}  

\theoremstyle{thmstyleone}%
\theoremstyle{thmstyletwo}%

\theoremstyle{thmstylethree}%

\begin{document}
	
	\title[Article Title]{A robust matrix-free approach for large-scale non-isothermal  high-contrast viscosity Stokes flow on blended domains with applications to geophysics}

	
	\author*[1]{\fnm{Andreas} \sur{Burkhart}\email{burk@cit.tum.de} \orcidlink{0009-0000-5460-5118}}
	
	\author[2]{\fnm{Nils} \sur{Kohl}\email{Nils.Kohl@lmu.de} \orcidlink{0000-0003-4797-0664}}
	
	\author[1]{\fnm{Barbara} \sur{Wohlmuth}\email{wohlmuth@cit.tum.de} \orcidlink{0000-0001-6908-6015}}
	
	\author[1]{\fnm{Jan} \sur{Zawallich}\email{jan.zawallich@cit.tum.de} \orcidlink{0000-0003-1602-2936}}
		
	\affil[1]{\orgdiv{Department of Mathematics, CIT}, \orgname{Technische Universit\"at M\"unchen}, \orgaddress{\city{Garching b. M\"unchen}, \postcode{85748}, \country{Germany}}}
	
	\affil[2]{\orgdiv{Department of Earth and Environmental Sciences}, \orgname{Ludwig-Maximilians-Universit\"at M\"unchen}, \orgaddress{\city{M\"unchen}, \postcode{80333}, \country{Germany}}}
	
	
	\abstract{
		We consider a compressible Stokes problem in the quasi-stationary case coupled with a time dependent advection-diffusion equation with special emphasis on high viscosity contrast geophysical mantle convection applications. In space, we use a P2-P1 Taylor--Hood element which is generated by a blending approach to account for the non-planar domain boundary without compromising the stencil data structure of uniformly refined elements. In time, we apply an operator splitting approach for the temperature equation combining the BDF2 method for diffusion and a particle method for advection, resulting in an overall second order scheme. Within each time step, a stationary Stokes problem with a high viscosity contrast has to be solved for which we propose a matrix-free, robust and scalable iterative solver based on Uzawa type block preconditioners, polynomial Chebyshev smoothers and a BFBT type Schur complement approximation. Our implementation is using a hybrid hierarchical grid approach allowing for massively parallel, high resolution Earth convection simulations.
	}
  	
	
	
	
	\maketitle
	
	\section{Introduction}\label{sec:Introduction}
	The Navier--Stokes equations coupled with an energy conservation equation \cite{Gassmoeller2020} can be used in geophysics to model the convective flow of the Earth's mantle as a viscous fluid. Since mantle convection takes place on a very long timescale ($\approx 10^8$ years), a common approach is to simplify the instationary Navier--Stokes equations to the stationary Stokes equations and argue that viscous relaxation and seismic waves effectively take place on a faster timescale than mantle convection \cite{Gassmoeller2020}. Further simplifications for a buoyancy driven fluid include the Boussinesq approximation \cite{Boussinesq1903} and the (truncated) anelastic liquid approximation \cite{Jarvis1980, Schubert2001}. 
    Such simplified Earth convection models can be classified into compressible \cite{Heister2017, Gassmoeller2020} or incompressible \cite{Kronbichler2012, Liu2019} approaches.	
	Many of the physical parameters driving the convective processes of a mantle convection simulation are hard to constrain and are thus only known up to a significant uncertainty. Prominent examples are the mantle viscosity for which estimates range from $10^{19}\uPa \us$ in the asthenosphere to $10^{23}\uPa \us$ in the lower mantle \cite[Chap.~5.7]{Fowler2004}, and the temperature at the core mantle boundary (CMB) for which estimates range from $2750\uK\pm 250\uK$ \cite{Jeanloz1986} to $4500\uK\pm 1000\uK$ \cite{daSilva2000, Jeanloz1990, Alfe1999}.
	    
	Of particular geophysical interest in terms of understanding geothermal events and the vertical motion of the Earth's lithosphere (uplift and subsidence) \cite{Bunge2018} is the study of rising hot material (called \textit{plumes}) and sinking cold material (called \textit{slabs}) inside the Earth's mantle. However, verifying the results of a given mantle convection model is a difficult task as we can only compare with observations based on the current state of the Earth's mantle, e.g., via seismic tomography which in turn introduces further uncertainty \cite{Freissler2024}. Comparing the result of a forward mantle convection simulation with real observations starting with our understanding of the mantle's current state as an initial state is also infeasible due the extremely long time scales at play.
	
	The aforementioned uncertainties and the fact that past states of the Earth's mantle can no longer be observed, motivates the study of inverse mantle convection problems using the adjoint method to facilitate reconstructions of the Earth’s mantle going backwards in time \cite{Colli2020}. Solving an inverse problem however typically includes a feedback loop using the result of forward model calculations to generate synthetic data and thus still requires fast solvers for the forward problem.
	
	The viscosity inside the Earth's mantle exhibits very steep gradients \cite[Chap.~5.7]{Fowler2004} and in general does not only depend on the location but also on the pressure, temperature and strain rate \cite{Gassmoeller2020}. This results in high contrast nonlinear partial differential equation (PDE) models and thus poses a significant numerical and mathematical challenge in particular to the iterative solvers typically used \cite{Kronbichler2012, Heister2017, Dannberg2022, Drzisga2018, Jodlbauer2024, Wichrowski2022}.
	
	In this work, we focus on a compressible forward mantle convection model using the truncated anelastic liquid approximation (TALA) coupled with a variation of the time dependent energy conservation equation presented in \cite{Gassmoeller2020}. We consider a space and temperature dependent nonlinear viscosity model. For the time dependent case, we use a semi-implicit time stepping scheme to linearize the resulting subproblems. In particular for the Stokes subproblem, we present an iterative flexible GMRES (FGMRES) \cite{Saad1993} outer Krylov solver preconditioned with an Uzawa type block preconditioner \cite{Drzisga2018}. The individual solver parts use geometric multigrid (GMG) approaches where applicable and are implemented in a matrix-free fashion in the finite element framework HyTeG\footnote{The open source software finite element framework HyTeG is freely available at \hyperlink{https://i10git.cs.fau.de/hyteg/hyteg}{https://i10git.cs.fau.de/hyteg/hyteg} under a GNU General Public License Version 3.} \cite{Kohl2019}, whose software architecture was designed for extreme parallel scalability \cite{Bauer2020}.	We show that this strategy can efficiently solve high resolution mantle convection problems despite the steep gradients in viscosity required for realistic geophysical simulations. 

    In the following, we introduce the problem setting (Sec.~\ref{sec:Model}), the time (Sec.~\ref{sec:Scheme}) and spatial discretisation (Sec.~\ref{sec:SpatialDisc}), the linear solvers and preconditioners (Sec.~\ref{sec:Variational} and \ref{sec:DiscTemp}). In Sec.~\ref{sec:GeodynamicApplication}, numerical tests are presented, including convergence and scalability tests, a coherent comparison of solvers and preconditioners, as well as high-resolution applications to real-world scenarios. This study is concluded in Sec.~\ref{sec:Conclusion}.
	\section{Model Problem} \label{sec:Model}
	We model the Earth's convective flow using the TALA \cite{Jarvis1980, Schubert2001, Gassmoeller2020} setting coupled to an energy conservation equation presented in \cite{Gassmoeller2020}.
    
    After a nondimensionalisation (see Appendix \ref{subsec:NondimEq}) similar to the one performed in \cite[Chap.~6.10]{Schubert2001} and \cite{Leng2008}, we arrive at
	\begin{align}
		& - \nabla \cdot \tau(u) + \nabla p_d  = -\frac{\Rayleigh}{\Peclet} g \rho\rk{p_s,T_s} \alpha\rk{p_s,T_s} T_d \, , \label{eq:momentum}
		\\
		& - \nabla \cdot u - \frac{\nabla \rho\rk{p,T}}{\rho\rk{p,T}} \cdot u = 0 \, , \label{eq:mass}
		\\
		& \rho\rk{p,T} C^p\rk{p,T} \left ( \frac{\partial T}{\partial t} + u \cdot \nabla T \right ) - \frac{1}{\Peclet} \nabla \cdot \left ( k\rk{p,T} \nabla T \right ) 
		\nonumber \\
		& - \Dissipation \alpha\rk{p,T} T \rho\rk{p_s,T_s} \left ( u \cdot g \right )  =  \rho\rk{p,T} H\rk{p,T} 
		+ \frac{\Peclet \Dissipation}{\Rayleigh} \tau(u) : \dot{\varepsilon}(u) \, , \label{eq:energy}
	\end{align}
	in $\left[0,T_{\text{end}} \right ] \times \Omega$, where $\Omega\subset\mathbb{R}^{\dim}$ is an approximation of the Earth's mantle, in our case
	\[
	\Omega=\{x\in\mathbb{R}^{\dim}\mid r_{\CMB} \leq \|x\|  \leq r_{\Surface}\},\quad r_{\CMB}<r_{\Surface},
	\]
	where we refer to the boundary as $\Gamma_{\CMB} \subset \partial \Omega$  if $\|x\|=r_{\CMB}$ and as $\Gamma_{\Surface} \subset \partial \Omega$ if $\|x\|=r_{\Surface}$. In case of the Earth's mantle, we have naturally $\dim=3$. However, we also consider $\dim=2$ for test cases. By an abuse of notation we denote the density as a function of $p_s, T_s$ and as a function of $p,T$ by the same symbol. We seek the velocity field $u\colon[0,T_{\text{end}}]\times\Omega\to\mathbb{R}^{\dim}$, the dynamic pressure $p_d\colon[0,T_{\text{end}}]\times\Omega\to\mathbb{R}$, and the temperature $T\colon[0,T_{\text{end}}]\times\Omega\to\mathbb{R}$. The viscous stress tensor $\tau$ and the deviatoric rate-of-deformation tensor $\dot{\varepsilon}$ are given by
	\[
		\tau(u):=2\eta(p,T,\dot{\varepsilon}(u))\dot{\varepsilon} \rk{u}, \quad \dot{\varepsilon}(u) := \nabla_S u -\frac{(\nabla\cdot u)I}{\dim} , \quad \nabla_S u := \frac{1}{2}\big(\nabla u+(\nabla u)^T\big),
	\]
	where $\eta$ is the viscosity, $\nabla_S u$ denotes the symmetric gradient, sometimes referred to as strain rate tensor \cite{Gassmoeller2020}.	 

    Starting from the coupled compressible Stokes equations presented in \cite[Sec.~2]{Gassmoeller2020}, we use a series of commonly used approximations to arrive at \eqref{eq:momentum}-\eqref{eq:energy} after nondimensionalisation. The TALA employs a first order Taylor approximation of the density in the buoyancy term of the momentum conservation equation to a given static reference pressure $p_s = p - p_d$ and temperature $T_s = T - T_d$. We use the hydrostatic approximation $\nabla p_s \approx \rho\rk{p_s, T_s}g$, see, e.g., \cite{King2010}, omit the density derivative in the mass conservation equation \cite[Sec.~1.2]{Gassmoeller2020}, omit the latent heating generated by phase transitions \cite[Sec.~2]{Gassmoeller2020} in the energy conservation equation and neglect in the part of the  adiabatic heating the dynamic pressure, i.e.,
    \[
		T \rk{ u \cdot \nabla p} \approx T \rk{ u \cdot \nabla p_s} \approx T \rho\rk{p_s, T_s} \rk{ u \cdot g}.
	\]
    The quantities
	$\eta\left ( p,T,\dot{\varepsilon}\left ( u \right ),x \right )$,
	$\rho\left ( p,T,x \right )$,	
	$k\left ( p,T,x \right )$,
	$H\left ( p,T,x \right )$,
	$\alpha\left ( p,T,x \right )$,
	$C^p\left ( p,T,x \right )$,
	$g\left ( x \right )$
	denote the viscosity, density, thermal conductivity, internal heating, thermal expansivity, specific heat capacity and gravity vector, respectively, and possibly show an explicit spatial dependency \cite{Gassmoeller2020}. Due to a simpler notation, we suppress the possible explicit spatial dependency in our notation in \eqref{eq:momentum}-\eqref{eq:energy} and in the following. The Rayleigh number $\Rayleigh$, Dissipation number $\Dissipation$ and Péclet number $\Peclet$ are dimensionless constants. 
	
	For simplicity, we assume a temperature and space dependent viscosity $\eta\colon \R \times \Omega \to \R$ and a space dependent density $\rho\colon\Omega\to\mathbb{R}$, neglecting the difference between $\rho\rk{p,T}$ and $\rho\rk{p_s,T_s}$.
    
    All scalar values $k$, $H$, $\alpha$, $C^p$ are assumed to be constant. The dimensional gravity vector is given by $-\frac{x}{\norm{x}} g_0$ for a reference constant $g_0 \in \R^+$ (usually $g_0 \approx 9.81$), hence in the nondimensional equations the gravity vector is given by $g\colon\Omega\to\R^{\dim}$ with $g\rk{x}\mapsto -\frac{x}{\norm{x}}$ (compare Appendix \ref{Nondim}). In a slight abuse of notation, we opted to use $p$ for the dynamic pressure $p_d$ in the following since our model no longer holds a reference to the total pressure. 

	Equations \eqref{eq:momentum} and \eqref{eq:mass} form a Stokes system driven by buoyancy and compressibility on the right-hand side. Equation \eqref{eq:energy} is a time dependent advection-diffusion equation driven by internal heating as the result of radioactive decay inside the Earth's interior, shear heating as the result of plastic deformation, and adiabatic heating as the result of compression or expansion of material. As boundary conditions for the velocity $u$, we prescribe Dirichlet boundary conditions at the surface boundary
	\[
	\evalat{u}_{\Gamma_{\Surface}} = u_{\Surface} \in \rk{ C^0\rk{\Gamma_{\Surface}}}^{\dim}
	\]
	and free slip boundary conditions at the CMB 
	\[
	\evalat{(u \cdot n)}_{\Gamma_{\CMB}} = 0 \text{ and } \evalat{{( \tau\rk{u} n)\cdot t\,}}_{\Gamma_{\CMB}} = 0,
	\]
	for all tangential vectors (on the boundary) $t$, 
	where $n$ is the outward facing unit normal vector. We note that $u_{\Surface}\cdot n=0$, and thus the boundary conditions are consistent with \eqref{eq:mass}, which can be reformulated to $\text{div}(\rho u)=0$. We enforce the first free-slip condition via a projection and the second free-slip condition is granted in a natural way by the weak formulation.
	The Dirichlet boundary condition at the surface can either be a no-slip boundary condition ($u_{\Surface}=0$) or can be chosen as tectonic plate velocities from tectonic plate reconstruction models. In case of using tectonic plate reconstructions, $u_{\Surface}$ is time dependent which we omit in our notation for simplicity. For our simulations, we are using the plate reconstruction model published in \cite{Mueller2022} using a scaling factor of $\frac{1}{2}$ to avoid forced convection \cite[Sec.~2.4]{Colli2020} ranging up to one billion years into the past. As boundary conditions for the temperature, we impose Dirichlet boundary conditions with a surface temperature of $T_{\Surface}=300$$\uK$ and a CMB temperature of $T_{\CMB}=4200$$\uK$ matching the pyrolite estimations of $4000\uK$--$4500\uK$ in \cite{daSilva2000}.
	
	\section{Time Discretisation and Semi-implicit Scheme} 
    \label{sec:Scheme}
	Let
	\begin{align*}
		& V :=  \gk{v \in \rk{H^1\rk{\Domain}}^{\dim} \text{ s.t. } \evalat{v}_{\Gamma_{\Surface}} = u_{\Surface} \text{ and } \evalat{\rk{v \cdot n}}_{\Gamma_{\CMB}} = 0 } ,
		\\ 
		& Q := \gk{ q \in L^2\rk{\Domain} \text{ s.t.} \int_\Domain q \dx = 0 }
        , 
        \\
        &  W := \gk{w \in H^1\rk{\Domain} \text{ s.t. } \evalat{w}_{\Gamma_{\Surface}} = \frac{T_{\Surface}}{\Delta T}\text{ and } \evalat{w}_{\Gamma_{\CMB}} = \frac{T_{\CMB}}{\Delta T} }
	\end{align*}
	denote the velocity, pressure and temperature function space, respectively. $\frac{T_{\Surface}}{\Delta T}$ and $\frac{T_{\CMB}}{\Delta T}$ denote the nondimensional surface and CMB temperature (compare Tab.~\ref{dimtable}). Discretising the time interval $[0,T_\text{end}]$ via $0=t^0<t^1<\ldots<t^N=T_\text{end}$ and setting $\tau^n=t^n-t^{n-1}$ for $n=1,\ldots,N$, let $u^n \in V,p^n \in Q$ and $T^n \in W$ denote the respective solutions at $t^n$.
    The Stokes system (equations \eqref{eq:momentum} and \eqref{eq:mass}) and the temperature equation \eqref{eq:energy} are solved in a Gauss--Seidel like alternating scheme. The energy equation is handled via a combined BDF2 and operator splitting approach.
    
    Similar to \cite{Kohl2022}, we introduce the characteristics $X: \ek{0, T_{\text{end}}}^2 \times \Omega \rightarrow \R^{\dim}$ with respect to a given velocity field $u\colon[0,T_{\text{end}}]\times\Omega\to\mathbb{R}^{\dim}$ and a fixed pair $\rk{s,x} \in \ek{0, T_{\text{end}}} \times \Omega$ as solutions to
    \begin{align*}
		\frac{\partial}{\partial t} X(s,t,x)  &= u\rk{X(s,t,x)}, \quad t \in \rk{0,T_{\text{end}}},\\
        X(s,s,x) &= x,
	\end{align*}
    where $X\rk{t_1,t_0,x}$ can be interpreted as the departure point we would arrive at if we start at $x$ and go backwards in time from $t^1$ to $t^0$ along $u$ \cite{Kohl2022}.
    
    For a fixed $s$, we can rewrite \eqref{eq:energy} in terms of
    \begin{align}
		\hat{T}\rk{t,x} := T\rk{t, X\rk{s,t,x}}
        , \quad
        \evalat{\frac{\partial}{\partial t}\hat{T}}_{\rk{s,x}} = \evalat{\rk{\frac{\partial}{\partial t} T + u \cdot \nabla T }}_{\rk{s,x}},\label{eq:timeDerivative}
    \end{align}
    resulting in
    \begin{align*}
		& \rho \rk{x} C^p \frac{\partial}{\partial t}\hat{T} - \frac{1}{\Peclet} \nabla \cdot \left ( k \nabla T \right ) 
	 - \Dissipation \alpha T \rho\rk{x} \left ( u \cdot g \right )  =  \rho\rk{x} H
		+ \frac{\Peclet \Dissipation}{\Rayleigh} \tau\rk{u} : \dot{\varepsilon}\rk{u} 
    \end{align*}
    for $t=s$ (compare \cite{Kohl2022}) together with the assumptions described in Sec.~\ref{sec:Model}.
    
    Setting $s = t^{n+1}$, we approximate  $\frac{\partial}{\partial t}\hat{T}$ using a variable time-step BDF2 \cite{Hairer1996} time discretisation, i.e.,
    \begin{align*}
		\frac{\partial}{\partial t}\hat{T} \approx \frac{D\ek{T^{n+1}, T\rk{t^n, X\rk{t^{n+1}, t^n, x}}, T\rk{t^{n-1}, X\rk{t^{n+1}, t^{n-1}, x}}}}{\tau^{n+1}}
    \end{align*}
    with 
    \begin{flalign*}
        D\ek{T^{n+1}, T^n, T^{n-1}} := \frac{2 \tau^{n+1} + \tau^{n}}{\tau^{n+1} + \tau^{n}} T^{n+1} - \frac{\tau^{n+1} + \tau^{n}}{\tau^{n}} T^{n} + \frac{\tau^{n+1} \cdot \tau^{n+1}}{\tau^{n} \left ( \tau^{n+1} + \tau^{n} \right )} T^{n-1}
    \end{flalign*}  
    (compare \cite[Sec.~3.1]{Kronbichler2012}). We refer to \cite{Etienne} and \cite{Becker} for stability considerations and constraints on $\frac{\tau^n}{\tau^{n+1}}$.
    
    The calculation of $\hat{T}\rk{t^n, x}$ and $\hat{T}\rk{t^{n-1}, x}$ is approximated via a particle based modified method of characteristics (MMOC) approach as introduced in \cite{Kohl2022}, compare \cite{Sime2021,Sime2022} for a different particle based approach in a similar setting. See \cite{Ettore} for an operator splitting in combination with streamlines tracking for a two-phase flow simulation in porous media. This allows us to handle advection-dominated systems without a stabilisation like the SUPG \cite{Brooks1982} or the entropy-viscosity method \cite{Guermond2011}.
    
    In particular, the MMOC assumes $u^{n+1}$, $u^n$ and $T^n$ to be known and tracks particles backwards in time along a linear interpolation between $u^{n+1}$ and $u^n$ using the classical Runge-Kutta method and evaluates $T^n$ at the calculated positions. Including the time step size $\tau^{n+1}$, we denote the application of the MMOC method by $\MMOC \rk{T^n, u^{n+1}, u^n, \tau^{n+1}}$. We approximate $\hat{T}\rk{t^n, x}$ and $\hat{T}\rk{t^{n-1}, x}$ by
    \begin{flalign*}
        \hat{T}\rk{t^n, x} &=T\rk{t^n, X\rk{t^{n+1}, t^n, x}}\approx \MMOC \rk{T^n, u_*^{n+1}, u^n, \tau^{n+1}},
        \\
        \hat{T}\rk{t^{n-1}, x} &\approx \MMOC \rk{ \MMOC \rk{T^{n-1}, u_*^{n+1}, u^n, \tau^{n+1}}, u^n, u^{n-1}, \tau^{n} },
    \end{flalign*}     
    while the rest of the splitting step, i.e., the energy equation without the advection, is solved via
    \begin{flalign}
        & \rho\left ( x \right ) C_p \left ( \frac{D\ek{T^{n+1}, \hat{T}^{n}, \hat{T}^{n-1}}}{\tau^{n+1}}  + { {u^{n+1}_* \cdot \nabla T^{n+1}}} \right )
        - { {\frac{1}{\text{Pe}}\nabla \cdot \left ( k \nabla T^{n+1} \right )}} 
        \nonumber \\ &
        - { {\text{Di} \, \alpha T^{n+1} \rho\left ( x \right ) \left ( u^{n+1}_* \cdot g \right ) }} 
        = { {{\rho\left ( x \right ) H}} }
        { +  {\frac{\text{Pe}\,\text{Di}}{\text{Ra}} 2\eta\left (x,T^{n+1}_* \right ) \dot{\varepsilon} \left ( u^{n+1}_* \right )  : {\dot{\varepsilon} \left ( u^{n+1}_* \right ) }}}, \label{eq:DiffSolve}
    \end{flalign}
    where
    \begin{flalign}
        & u^{n+1}_* := u^{n} + \frac{u^{n} - u^{n-1}}{\tau^{n}} \tau^{n+1}, \quad T^{n+1}_* := T^{n} + \frac{T^{n} - T^{n-1}}{\tau^{n}} \tau^{n+1}
        \label{eq:Extrapolation}
    \end{flalign}
    are linear extrapolations in time. With the temperature at the new time step, we solve the Stokes system via
    \begin{equation}\label{eq:StokesSub}
    \begin{split}
        - \nabla \cdot 2\eta(x,T^{n+1})\dot{\varepsilon}(u^{n+1}) + \nabla p^{n+1}  &= { -\frac{\text{Ra}}{\text{Pe}} g  \rho \left ( x \right ) \alpha {T_d^{n+1}} }, \\ - \nabla \cdot u^{n+1} -  \frac{\nabla \rho(x)}{\rho(x)} \cdot u^{n+1} &= 0.
    \end{split}
    \end{equation}	
    To generate a first step in the multi-step method \eqref{eq:DiffSolve}, we use the implicit Euler method. We assume $u^0$ to be zero. Therefore, the Stokes system in this step is incompressible. 
    In each iteration, a new time step $\tau^{n+1}$ is selected, satisfying the previously mentioned stability constraints as well as a Courant--Friedrichs--Lewy condition \cite{Quarteroni2008}, s.t.
	\begin{equation*}
		\tau^{n+1}\max_{K\in\mathcal{T}_{h}}\frac{\|u^n\|_{L^\infty(K)}}{h_K} \leq C_{CFL}
	\end{equation*} 
	where $\|u^n\|_{L^\infty(K)}$ denotes the maximal magnitude of the velocity on the element $K$ at time $t^n$, and $h_K$ denotes the element diameter \cite{Kronbichler2012}, both introduced in the following section. The parameter $C_{CFL}$ is chosen experimentally. Note that in practice, the curvature of the domain (compare Sec.~\ref{subsec:Mesh}) can be neglected when calculating $h_K$.
	\section{Spatial Discretisation} \label{sec:SpatialDisc}
    In the following, we describe the spatial discretisation, including how the Earth's mantle is approximated via a mesh and blending, the finite element method, and the resulting discrete variatonal problems.
	\subsection{Mesh and Blending} \label{subsec:Mesh}
	We discretise $\Omega$ using a hybrid hierarchical grid approach native to the HHG \cite{Bergen2004, Bergen2006, Gmeiner2012,Gmeiner2014,Gmeiner2015,Gradl2015,Gmeiner2016} and HyTeG \cite{Kohl2019} finite element frameworks. This hybrid discretisation approach uses an unstructured coarse grid consisting of triangles in 2D and tetrahedrons in 3D on which a structured refinement is performed \cite{Bey1995}.
	
	In case of the annulus, we start with the generation of an unstructured grid by dividing a hollow convex regular polygon into congruent isosceles trapezoids. The trapezoids are then divided into layers of equal height. The trapezoids in these layers are then in turn subdivided into triangles. The generation of the unstructured hollow spherical mesh works in a similar way. We start with a triangular truncated pyramid mesh, based on the hollow icosahedron. These truncated pyramids are again divided into layers of equal height into congruent truncated pyramids, which are then in turn subdivided into tetrahedrons in each layer.
	
\begin{figure}[htbp]
		\centering
		\includegraphics[width=0.9\textwidth]{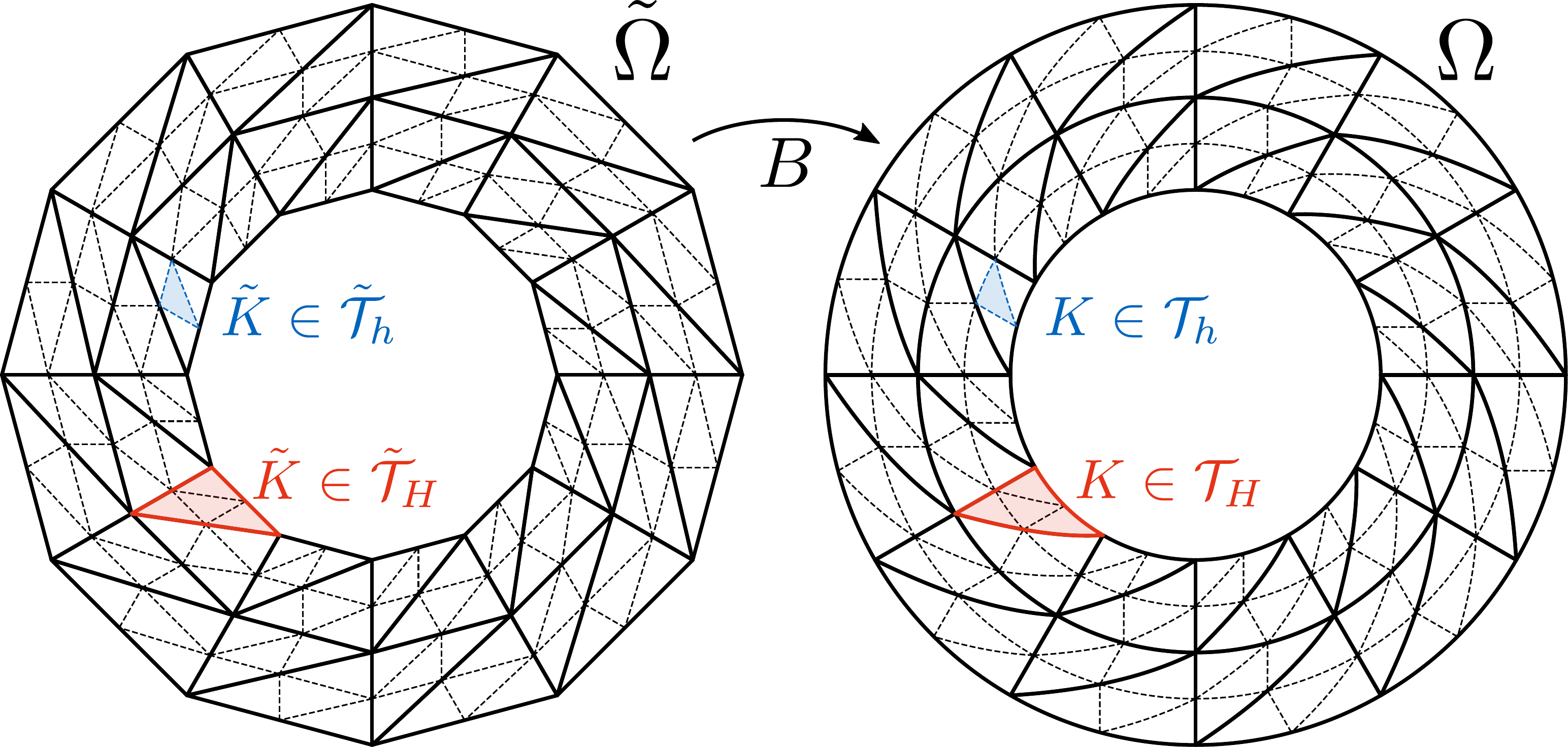}
		\caption{Mapping of $\tilde{\Omega}$ covered by a triangular mesh with a structured refinement (left) onto an annulus $\Omega$ (right) via an appropriate blending map. The first level of grid refinement is indicated with dashed lines. $\tilde{K}$ denotes an exemplary element on the unstructured coarse grid $\tilde{\mathcal{T}}_{H}$ and refined grid $\tilde{\mathcal{T}}_{h}$. The respective blending map images are denoted as $K$ (compare Fig.~\ref{fig:BlendingTransformation}).}
        \label{fig:BlendedDomain}
	\end{figure}

	In order to accurately represent the geometry of the refined unstructured mesh, it is mapped to the physical domain $\Omega$ by a blending map as depicted in Fig.~\ref{fig:BlendedDomain}. Let $\tilde{\Omega} := \bigcup_{\tilde{K} \in \tilde{\mathcal{T}}_{h}} \tilde{K}$ be the set covered by the elements $\tilde{K}$ of the refined unstructured mesh $\tilde{\mathcal{T}}_{h}$. Let the blending map $B\colon\tilde{\Omega} \rightarrow \Omega$ be a $C^0$-diffeomorphism (which is locally a $C^1$-diffeomorphism on each macroelement) and let $\mathcal{T}_{h} := \gk{B\rk{\tilde{K}}\mid \tilde{K} \in \tilde{\mathcal{T}}_{h}}$ denote the set of blended elements $K \in \mathcal{T}_{h}$.
    
    \begin{figure}[htbp]
		\centering
		\includegraphics[width=0.9\textwidth]{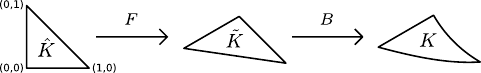}
		\caption{Mapping of the 2D reference triangle element $\hat{K}$ to a blended element $K$ in the domain via the application of an affine transformation $F$ and a blending transformation $B$.} \label{fig:BlendingTransformation}
	\end{figure} 	

    Fig.~\ref{fig:BlendingTransformation} shows the process of mapping a reference element $\hat{K}$ onto a blended element $K$ with an intermediate affine transformation step. The evaluation of integrals over a blended element $K$ pulled back to a reference element $\hat{K}$ as typically performed in a finite element context needs to take the blending map into account and requires the evaluation of the blending Jacobian $J_B$. For example, a prototypical diffusion integral in the discrete velocity space takes the form
    \begin{flalign*}
		&\int_{K} \nabla_x \varphi_i^u : \nabla_x \varphi_j^u \dx =

        \int_{\hat{K}}
        
        \nabla_x
          \varphi_i^{\hat{u}}
        
        :
        
        \nabla_x
          \varphi_i^{\hat{u}} 
        
        \Big| \det (\evalat{J_B}_{F \rk{\hat{x}}}) \Big|
         \Big|  \det (\evalat{J_F}_{\hat{x}}) \Big| \dc{\hat{x}}
	\end{flalign*}
    with 
    \begin{flalign*}
		\nabla_x \varphi_i^{\hat{u}} =\rk{
        \evalat{J_B^{-1}}_{F\rk{\hat{x}}}
        }^{T} 
        \rk{\evalat{J_F^{-1}}_{\hat{x}}}^{T} \nabla_{\hat{x}}\varphi_i^{\hat{u}}
    \end{flalign*}    
    and $\varphi_i^{\hat{u}}$ representing the nodal basis functions associated with the reference element. 

	\subsection{Finite Element Method} \label{subsec:FEM}
	We spatially approximate the solution functions on the refined unstructured mesh $\tilde{\mathcal{T}}_{h}$ using the finite element method. For 
    $\rk{u, p}$ we use $\mathbb{P}^{\dim}_2 \times \mathbb{P}_1$ Taylor--Hood elements and for $T$ conforming $\mathbb{P}_2$ elements. In particular, we are using nodal basis functions associated with the Lagrange interpolation points of order 1 and 2. In the following we will denote the vertex and edge midpoint positions defined by the Lagrange interpolation points of order 2 as 
    \begin{equation}\label{eq:Lagrange}
        \tilde{x}_i \in \tilde{\Domain}, \, i \in I_{\tilde{x}} := I_{\Inner} \dot{\cup} I_{\Surface} \dot{\cup} I_{\CMB}, 
    \end{equation}
    which we split into inner, surface and CMB positions.
    
    The Taylor--Hood element pairing as introduced in \cite{Taylor1973}, is known to be a discrete inf-sup stable pairing for the Stokes problem under certain conditions on the domain and mesh \cite{Verfuerth1984, Girault1986, Ern2004} and is one of the default element choices for geophysical applications \cite{Thieulot2022}. 
	\section{Discrete Stokes Problem and Solver} \label{sec:Variational}
	Taking into account blending, we define the discrete function spaces
    \begin{flalign*}
    V_h = \Big\{ & v_h \mid v_h = \tilde{v}_h \circ B^{-1} \text{ for some }
    \\
    & \tilde{v}_h \in \rk{C^0\rk{\tilde{\Domain}}}^{\dim} \text{ such that }  
    \forall \tilde{K} \in \tilde{\mathcal{T}}_h \, \evalat{\tilde{v}_h}_{\tilde{K}} \in \rk{\mathcal{P}_2\rk{\tilde{K}}}^{\dim}  \Big\},
    \\
    Q_h = \Big\{ & q_h \mid q_h = \tilde{q}_h \circ B^{-1} \text{ for some }
    \\
    & \tilde{q}_h \in L^2\rk{\tilde{\Domain}} \cap C^0\rk{\tilde{\Domain}} \text{ such that } 
    \forall \tilde{K} \in \tilde{\mathcal{T}}_h \, \evalat{\tilde{q}_h}_{\tilde{K}} \in \mathcal{P}_1\rk{\tilde{K}}  \Big\},
    \\
    W_h = \Big\{ & w_h \mid w_h = \tilde{w}_h \circ B^{-1} \text{ for some }
    \\
    & \tilde{w}_h \in C^0\rk{\tilde{\Domain}} \text{ such that }  \forall \tilde{K} \in \tilde{\mathcal{T}}_h \, 
    \evalat{\tilde{w}_h}_{\tilde{K}} \in \mathcal{P}_2\rk{\tilde{K}}  \Big\} ,
    \end{flalign*}
    with $\mathcal{P}_k \left ( \tilde{K} \right )$ representing the space of multivariate polynomials up to global degree $k$ on an element $\tilde{K}\in\tilde{\mathcal{T}}_h$.
    \subsection{Linear algebraic formulation} 
    We seek to find approximate solutions $u_h^{n+1} = \tilde{u}_h^{n+1} \circ B^{-1} \in V_h$, $p_h^{n+1} \in Q_h$ to \eqref{eq:StokesSub} that satisfy the boundary and pressure constraints 
    \begin{gather*} 
    \forall i \in I_{\CMB} \, \evalat{\rk{u_h^{n+1} \cdot n}}_{B\rk{\tilde{x}_i}} = 0, \, \, \, \forall i \in I_{\Surface}\, \tilde{u}^{n+1}_h\rk{\tilde{x}_i} = u_{\Surface}\rk{B\rk{\tilde{x}_i}}, \, \, \,
    \int_{\Domain} p_h^{n+1} \dc{x} = 0.
    \end{gather*}
    Let $\varphi_{0}^u, \ldots, \varphi_{N_u-1}^u$ and $\varphi_{0}^p,\ldots,\varphi_{N_p-1}^p$ denote the velocity and pressure nodal basis functions associated with the blended Taylor--Hood element on $\Domain$.
    
    Neglecting the boundary conditions at first, the discrete variational Stokes problem can be written as a linear system using ${\textbf{f}}^{\,n+1}\in \R^{N_u}$, ${\textbf{g}}=\mathbf{0} \in \R^{N_p}$, ${B},{C} \in \mathbb{R}^{N_p \times N_u}$, and symmetric, positive semidefinite ${A} \in \mathbb{R}^{N_u \times N_u}$, where    
	\begin{flalign*}
		{A}_{ij} = &\int_\Omega 2 \eta\rk{x,T_h^{n+1}} \nabla_s \varphi_j^u : \nabla_s \varphi_i^u \dx
		- \frac{2}{\dim} \int_\Omega \rk{\nabla \cdot \varphi_j^u}\rk{\nabla \cdot \varphi_i^u} \dx,
		\\
		{B}_{ij} =&  - \int_\Omega \rk{\nabla  \cdot \varphi_j^u} \varphi_i^p\dx, \quad {C}_{ij} =  - \int_\Omega \nabla \ln \rho \cdot \varphi_j^u \varphi_i^p \dx,
		\\
		{\textbf{f}}^{\,n+1}_i = & - \int_\Omega \frac{\text{Ra}}{\text{Pe}} \rho \alpha {T_{d,h}^{n+1}} g \cdot \varphi_i^u \dx,
	\end{flalign*}
    resulting in an asymmetrical generalised saddle point problem with block structure
    	\begin{flalign}
    		\begin{pmatrix} {A} & {B}^T \\ {B}+{C} & 0 \end{pmatrix} \begin{pmatrix} \,\,{\textbf{u}}^{n+1}\,\,\\ \,\, {\textbf{p}}^{n+1}\,\, \end{pmatrix} = \begin{pmatrix} \,\,{\textbf{f}}^{\,n+1} \,\,\,\\ \, \, {\textbf{g}} \,\,\, \end{pmatrix} \label{eq:SaddlePointStructure}
    	\end{flalign}
    and coefficient vectors ${\textbf{u}}^{n+1} \in \R^{N_u}$, ${\textbf{p}}^{n+1} \in \R^{N_p}$. 

    As a second step, we assume the velocity degrees of freedom to be ordered such that, similar to \eqref{eq:Lagrange}, we can split up 
    \begin{flalign}
        \mathbf{u}^{n+1} =: \mat{
        \mathbf{u}_{\Surface}^{n+1}
        \\[2pt]
        \mathbf{u}_{\CMB}^{n+1}
        \\[2pt]
        \mathbf{u}_{\Inner}^{n+1}
        }
        ,\quad
        {\textbf{f}}^{\,n+1} =: \mat{
        {\textbf{f}}_{\Surface}^{\,n+1}
        \\[2pt]
        {\textbf{f}}_{\CMB}^{\,n+1}
        \\[2pt]
        {\textbf{f}}_{\Inner}^{\,n+1}
        },
    \end{flalign}
    between surface, CMB and inner degrees with $\mathbf{u}_{\Surface}^{n+1}, {\textbf{f}}_{\Surface}^{\,n+1} \in \R^{\dim\abs{I_{\Surface}}}$, $\mathbf{u}_{\CMB}^{n+1}, {\textbf{f}}_{\CMB}^{\,n+1} \in \R^{\dim\abs{I_{\CMB}}}$ and $\mathbf{u}_{\Inner}^{n+1}, {\textbf{f}}_{\Inner}^{\,n+1} \in \R^{\dim\abs{I_{\Inner}}}$, respectively. Splitting up 
    \begin{flalign*}
		\mat{
        {A}_{SS} & {A}_{SC} & {A}_{SI}  \\
        {A}_{CS} & {A}_{CC} & {A}_{CI}  \\
        {A}_{IS} & {A}_{IC} & {A}_{II}  \\
        } := {A} , \quad 

        \mat{{B}_{IS} & {B}_{IC} & {B}_{II}} := {B} , \quad

        \mat{{C}_{IS} & {C}_{IC} & {C}_{II}} := {C},
	\end{flalign*}
    into blocks of matching dimensions results in a block system matrix of the form
    \begin{flalign}
        {M} := \mat{
        {A}_{SS} & {A}_{SC} & {A}_{SI} &  {B}_{IS}^T   \\[2pt]
        {A}_{CS} & {A}_{CC} & {A}_{CI} &  {B}_{IC}^T   \\[2pt]
        {A}_{IS} & {A}_{IC} & {A}_{II} &  {B}_{II}^T \\[2pt]
        \bar{B}_{IS} & \bar{B}_{IC} & \bar{B}_{II} & 0
        } \in \R^{\rk{N_u + N_p} \times \rk{N_u + N_p}}.
    \end{flalign}
    with $\bar{B}:=B+C \in \mathbb{R}^{N_p \times N_u}$.
    
    The pressure solution is only determined up to a constant (compare \cite{Ern2004, Boffi2008}) and the surface normal on the CMB in general does not align with the coordinate axes. Hence the free-slip boundary condition cannot be enforced individually in each dimension and needs to be treated in a different fashion than the surface boundary condition. To address this, as a final step, we introduce a discrete projection operator
    \begin{flalign}
        P := \mat{
        I_{S} & 0   & 0   & 0   \\
        0   & P_v & 0   & 0   \\
        0   & 0   & I_{I} & 0   \\
        0   & 0   & 0   & P_p
        } \in \R^{\rk{N_u + N_p} \times \rk{N_u + N_p}}
    \end{flalign}
    with identity matrices $I_{S} \in \R^{\rk{\dim\abs{I_{\Surface}}} \times \rk{\dim\abs{I_{\Surface}}}}$, $I_{I} \in \R^{\rk{\dim\abs{I_{\Inner}}} \times \rk{\dim\abs{I_{\Inner}}}}$ and projections $P_p \in \R^{N_p \times N_p}$, $P_v \in \R^{\rk{\dim\abs{I_{\CMB}}} \times \rk{\dim\abs{I_{\CMB}}}}$, respectively ensuring a unique pressure solution and that the free-slip boundary condition at the CMB is discretely fulfilled. Note that in practice $P_p$ is implemented such that it only enforces that the vector $\textbf{p}^{n+1}$ has arithmetic mean zero. This saves computational effort during the solving process and an additive constant $c_p \in \R$ to ensure $\int_{\Domain} \rk{p_h^{n+1} + c_p} \dx = 0$ can be determined in a post processing step. 
    
    Applying $P$ to ${M}$ in a similar fashion to a left and right preconditioner, whilst realising the Dirichlet boundary conditions on $\Gamma_{\Surface}$ by a standard block row modification results in a system of the form
    \begin{flalign}
        P \mat{
        I_{S} & 0 & 0 &  0   \\[2pt]
        {A}_{CS} & {A}_{CC} & {A}_{CI} &  {B}_{IC}^T   \\[2pt]
        {A}_{IS} & {A}_{IC} & {A}_{II} &  {B}_{II}^T \\[2pt]
        \bar{B}_{IS} & \bar{B}_{IC} & \bar{B}_{II} & 0
        } P \mat{
        \mathbf{u}_{\Surface}^{n+1}
        \\[2pt]
        \mathbf{u}_{\CMB}^{n+1}
        \\[2pt]
        \mathbf{u}_{\Inner}^{n+1}
        \\[2pt]
        \mathbf{p}^{n+1}
        } = P \mat{
        {\textbf{u}}_{\Surface}^{\text{Int}}
        \\[2pt]
        {\textbf{f}}_{\CMB}^{\,n+1}
        \\[2pt]
        {\textbf{f}}_{\Inner}^{\,n+1}
        \\[2pt]
        \mathbf{g}
        } \label{eq:SaddlePointStructureDirichlet}
    \end{flalign}
    where ${\textbf{u}}_{\Surface}^{\text{Int}}$ denotes the interpolation $I_h \rk{u_{\Surface}}$ of the surface boundary condition.  Note that $\textbf{u}^{n+1}$ and $\textbf{p}^{n+1}$ in \eqref{eq:SaddlePointStructureDirichlet} are only determined up to the projections $P_v$ and $P_p$, respectively, hence \eqref{eq:SaddlePointStructureDirichlet} is singular. Using an FGMRES solver to construct a Krylov subspace basis \cite[Sec.~2.1]{Saad1993} and starting with a consistent (projected) initial guess and right-hand side, this problem is circumvented, leading the approximated solution to automatically fulfill the discrete free-slip boundary condition and to have pressures of mean zero. This approach however requires careful application of the pressure and free-slip projections throughout the solver implementation. In particular, we apply these projections as part of the block preconditioner (see Sec.~\ref{subsec:UzawaPreconditioner}), multigrid prolongation, restriction and smoother (see Sec.~\ref{subsec:AApproximation}). 
    
	\subsection{Stokes System Solver} \label{sec:StokesSolvers}
    Throughout our model implementation in HyTeG, we use state-of-the-art automated code generation introduced in \cite{Boehm2025} to generate vectorized C++ code using optimized elementwise looping strategies in order to efficiently evaluate the application of linear operators\footnote{We use a modified version of the HyTeG operator generator, see \href{https://i10git.cs.fau.de/hyteg/hog}{https://i10git.cs.fau.de/hyteg/hog}, that automatically creates vectorized C++ code for the application of linear operators.}.

    In this section we focus on the generalised saddle point sub system of \eqref{eq:SaddlePointStructureDirichlet} with respect to the inner and CMB degrees of freedom, which we rewrite as
    \begin{flalign}
        \mat{
        P_v {A}_{CC} P_v & P_v {A}_{CI} &  P_v {B}_{IC}^T P_p  \\[2pt]
        {A}_{IC} P_v  & {A}_{II} &  {B}_{II}^T P_p \\[2pt]
        P_p \bar{B}_{IC} P_v  & P_p \bar{B}_{II} & 0
        } \mat{
        \mathbf{u}_{\CMB}^{n+1}
        \\[2pt]
        \mathbf{u}_{\Inner}^{n+1}
        \\[2pt]
        \mathbf{p}^{n+1}
        } = \mat{
        P_v \rk{{\textbf{f}}_{\CMB}^{\,n+1} - {A}_{CS} {\textbf{u}}_{\Surface}^{\text{Int}}}
        \\[2pt]
        {\textbf{f}}_{\Inner}^{\,n+1} - {A}_{IS} {\textbf{u}}_{\Surface}^{\text{Int}}
        \\[2pt]
        P_p \rk{\mathbf{g} - \bar{B}_{IS} {\textbf{u}}_{\Surface}^{\text{Int}}}
        }. \label{eq:InnerSaddlePointSystem}
    \end{flalign}
    In an abuse of notation with respect to \eqref{eq:SaddlePointStructure} we will denote \eqref{eq:InnerSaddlePointSystem} as
    \begin{flalign}
        \mat{
        A & B^T \\[2pt]
        B+C & 0
        } \mat{
        \mathbf{u}^{n+1}
        \\[2pt]
        \mathbf{p}^{n+1}
        } = \mat{{\textbf{f}}^{\,n+1}
        \\[2pt]
        \mathbf{g}
        } \label{eq:InnerSaddlePointSystemGeneralised}
    \end{flalign}
    with $A \in \mathbb{R}^{\rk{N_u - \dim\abs{I_{\Surface}}} \times \rk{N_u - \dim\abs{I_{\Surface}}}}$, ${B},{C} \in \mathbb{R}^{N_p \times \rk{N_u - \dim\abs{I_{\Surface}}}}$, $\mathbf{u}^{n+1}, {\textbf{f}}^{\,n+1} \in \R^{\rk{N_u - \dim\abs{I_{\Surface}}}}$ and $\mathbf{p}^{n+1}, \mathbf{g} \in \R^{N_p}$ in the following. In this context, inverse operators like $A^{-1}$ or the inverse Schur complement $S^{-1}$ (see Sec.~\ref{subsec:UzawaPreconditioner}) as well as their respective approximations are to be interpreted as well-defined operations, e.g., via the application of projections and Krylov solvers starting with a consistent initial residual.
	\subsubsection{Uzawa Type Block Preconditioners} \label{subsec:UzawaPreconditioner}
	We make use of Uzawa type block preconditioners \cite{Drzisga2018} for the compressible Stokes system. In particular, we consider the inexact Uzawa, adjoint inexact Uzawa, and symmetric Uzawa block preconditioners \cite[Sec.~3]{Drzisga2018}, as defined for a symmetrical saddle point system, i.e., for $C=0$. Omitting the time step for a clearer notation, a single update step of these block preconditioners for the saddle point system can be written as \cite[Sec.~3.1]{Drzisga2018}:
	\\ \\
	Inexact Uzawa:
	\begin{equation}
    \begin{split}
		& \textbf{u}_{k+1} := \textbf{u}_k + \sigma \hat{A}^{-1} \left ( \textbf{f} - A \textbf{u}_k - B^T \textbf{p}_{k} \right )
		\\
		& \textbf{p}_{k+1} = \textbf{p}_k - \omega \hat{S}^{-1} \left ( \textbf{g} - B \textbf{u}_{k+1} \right )
        \end{split}
        \label{eq:InexactUzawaStep1}
	\end{equation}	
	Adjoint inexact Uzawa:
	\begin{equation}
    \begin{split}
		& \textbf{p}_{k+1} = \textbf{p}_k - \omega \hat{S}^{-1} \left ( \textbf{g} - B \textbf{u}_k \right ) 
		\\
		& \textbf{u}_{k+1} := \textbf{u}_k + \sigma \hat{A}^{-1} \left ( \textbf{f} - A \textbf{u}_k - B^T \textbf{p}_{k+1} \right )
    \end{split}\label{eq:AdjointInexactUzawaStep1}
	\end{equation}	
	Symmetric Uzawa:
	\begin{equation}
    \begin{split}
		& \hat{\textbf{u}} := \textbf{u}_k + \sigma \hat{A}^{-1} \left ( \textbf{f} - A \textbf{u}_k - B^T \textbf{p}_k \right ) 
		\\
		& \textbf{p}_{k+1} = \textbf{p}_k - \omega \hat{S}^{-1} \left ( \textbf{g} - B \hat{\textbf{u}} \right ) 
		\\
		& \textbf{u}_{k+1} := \hat{\textbf{u}} + \sigma \hat{A}^{-1} \left ( \textbf{f} - A \hat{\textbf{u}} - B^T \textbf{p}_{k+1} \right )
    \end{split}
    \label{eq:SymmetricUzawaStep1}
	\end{equation}	
	Here, $\hat{A}$ and $\hat{S}$ are representing nonsingular approximations of $A$ and the Schur complement $S = B A^{-1} B^T$ with respect to the symmetrical system. 
    The Uzawa method can also be used as a saddle point system solver \cite{Elman1994} and smoother \cite{Drzisga2018}. The adjoint inexact Uzawa matches the block triangular preconditioner used in \cite{Kronbichler2012}. Regarding the convergence of the block preconditioner, the relaxation parameters $\omega, \sigma > 0$ should be chosen appropriately such that $\sigma \hat{A} \geq A$ and $\omega \hat{S} \geq S$ \cite{Drzisga2018, Jodlbauer2024}. Note that in the block preconditioners above, $B$ can be replaced by $B+C$ while not replacing $B^T$ by $(B+C)^T$, resulting in a non-symmetrical approach with a slightly higher convergence rate for some examples.
	\subsubsection{Approximation of $A^{-1}$} \label{subsec:AApproximation}
	In order to approximate $A^{-1}$, we use a conjugate gradient (CG) solver preconditioned by a GMG V-cycle solving the $A$ block up to a relative residual tolerance $\tol_A>0$. We are using $m_A\in\mathbb{N}$ pre- and postsmoothing steps of a Chebyshev smoother \cite{Jodlbauer2024} of polynomial degree $\degree_A\in\mathbb{N}$ preconditioned via a multiplication with the inverse diagonal of $A$ as a smoother for the multigrid approach on each refinement level. Note that in order to avoid division by zero when the surface normal on the CMB aligns with the coordinate axes we exclude the free-slip projection $P_v$ during the calculation of $\diag\rk{A}$.
    
    On the coarsest grid, we use a non-matrix free CG solver preconditioned by an algebraic multigrid (AMG) method implemented in the freely available software package PETSc\footnote{The open source software package PETSc for solving coarse grid problems is freely available at \hyperlink{https://petsc.org}{https://petsc.org} under a 2-clause BSD license.} solving up to a relative tolerance of $\tol_{\coarse}>0$. As an alternative to PETSc, a CG solver with optional agglomeration of the coarse grid problem to a lower number of processes implemented in HyTeG \cite{Huber2020} can be used. 

    \begin{figure}[htbp]
    \centering
    \includegraphics[width=0.55\textwidth]{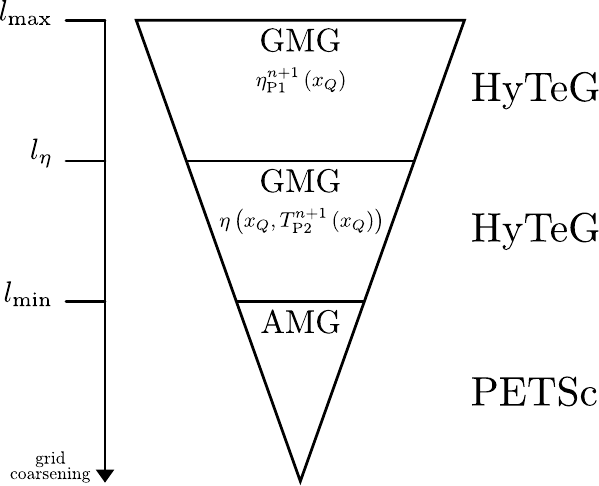}
		\caption{Solver structure of the $A$ block V-cycle used to precondition the respective CG solver. $l_{\max}$ and $l_{\min}$ denote the finest and coarsest grid refinement level, respectively. On levels greater than $l_{\min}$, we use a GMG approach, while an AMG approach is used on the coarse grid.}
		\label{fig:ASolverStructure}
	\end{figure}	
    Since the evaluation of the exponential function (see choice of viscosity in Sec.~\ref{subsec:GeoModel}) is difficult to vectorize, we have opted to represent the viscosity as a P1 function. On coarse refinement levels, however, a P1 function might not be able to resolve high contrasts in the viscosity accurately, negatively impacting the convergence rate of the solver. On refinement levels up to level $l_{\eta}$ (with lower levels being coarser), we first calculate an evaluation $T^{n+1}_{\text{P2}}\rk{x_Q}$ of the temperature function at quadrature points $x_Q \in \Domain$ and calculate $\eta\rk{x_Q, T^{n+1}_{\text{P2}}\rk{x_Q}}$ as the evaluation of the viscosity at $x_Q$. On refinement levels above $l_{\eta}$, we precalculate the viscosity $\eta^{n+1}_{\text{P1}}$ as a P1 function with values $\eta\rk{x_V, T^{n+1}_{P2}\rk{x_V}}$ at the vertex locations $x_V \in \Domain$ and evaluate $\eta^{n+1}_{\text{P1}}$ at the quadrature points as part of applying the operator $A$. Fig.~\ref{fig:ASolverStructure} depicts the overall structure of the V-cycle used to precondition the $A$ block CG solver with respect to the maximum and minimum grid refinement level $l_{\max}$, $l_{\min}$.

    Since the precalculation of $\eta^{n+1}_{\text{P1}}$ needs to be performed only once per time step, a vectorized evaluation of the exponential function on levels $l_{\eta}+1, \ldots, l_{\max}-1, l_{\max}$ is non-essential. In order to vectorize the evaluation of the exponential function on levels $l_{\eta}$ and below, we use an operator generated and optimized for the specific viscosity \eqref{eq:ViscosityDef} with $\eta_{\base}$ as depicted in Fig.~\ref{fig:baseViscosity}. For this operator, we replace the exponential function with a piecewise polynomial approximation of order 6 with an relative approximation error smaller than $6\cdot10^{-4}$ in the required range of evaluation. The errors introduced by this approximation are smoothed out on levels above $\l_{\eta}$ and do not negatively impact the convergence rate of the solver.
    
    For arbitrary choices of $\eta\rk{x,T}$ that can not be well approximated by polynomials, we provide in our code a less efficient and scalable $A$ block operator implementation on levels $l_{\eta}$ and below that we by default apply in a non-matrix free fashion.
	\subsubsection{Approximation of the Inverse Schur Complement} \label{subsec:SchurApproximation}
	A common choice, see \cite{Burstedde2008, Geenen2009, Burstedde2009, Furuichi2011, Kronbichler2012, Burstedde2013, May2015, Isaac2015, Gmeiner2016}, for the Schur complement approximation $\hat{S}$ in case of sufficiently smooth viscosities, is the inverse-viscosity scaled mass matrix $M_{1/\eta}$ (or a lumped version of it), i.e.,
	\begin{flalign}\label{eq:schurMass}
		\hat{S}_M^{-1} :=  M^{-1}_{1/\eta}.
	\end{flalign}
	The mass approximation was shown to be spectrally equivalent to the Schur complement in the isoviscous \cite{Verfuerth1984b} and variable viscosity case \cite{Grinevich2009}, but is known to deliver poor convergence results in case of high viscosity contrasts \cite{Rudi2015, Rudi2017}. In our tests, we solve $\hat{S}_M$ up to a relative tolerance of $\tol_{\invMass}>0$.

	Another choice are BFBT (or least-squares commutator) type inverse Schur complement approximations \cite{Elman1999, Elman2006a, Elman2006, May2008, Furuichi2011,Isaac2015,Rudi2015,Rudi2017} of the form
	\begin{flalign*}
		(BA^{-1}B^T)^{-1}=\left ( B {C}^{-1} B^T \right )^{-1} \left ( B {C}^{-1} A {D}^{-1} B^T  \right ) \left ( B {D}^{-1} B^T \right )^{-1}
	\end{flalign*}
	for suitable symmetric positive definite matrices $C,D$. Prominent examples are the diag(A)-BFBT approximation, where $C=D=\diag(A)$, and the weighted BFBT approximation
	\begin{flalign} \label{eq:schurWBFBT}
		\hat{S}_{w}^{-1} := \left ( B {M}_{\sqrt{\eta}}^{-1} B^T \right )^{-1} \left ( B {M}_{\sqrt{\eta}}^{-1} A {M}_{\sqrt{\eta}}^{-1} B^T  \right ) \left ( B {M}_{\sqrt{\eta}}^{-1} B^T \right )^{-1}
	\end{flalign}
	with ${M}_{\sqrt{\eta}}$ as the (potentially lumped) vector mass scaled with the square root of the viscosity, which have been shown to be a  suitable and algorithmically robust choice, especially for geodynamic applications with high viscosity contrasts \cite{Rudi2015, Rudi2017}. In particular, the weighted BFBT approximation has been shown to preserve its efficacy in combination with higher order space discretisations using an asymmetric rescaling of the viscosity close to Dirichlet boundaries for mesh independent algorithmic robustness \cite{Rudi2017}. Similar to \cite{Rudi2015, Rudi2017}, we replace $\left ( B {M}_{\sqrt{\eta}}^{-1} B^T \right )^{-1}$ by solving an inverse viscosity scaled Poisson problem with Neumann boundary conditions up to a relative or absolute tolerance of $\tol_{\wBFBT}>0$ (whichever applies first) with a CG solver preconditioned by a GMG approach which we denote by $K^{-1}_{1/\sqrt{\eta}}$ and introduce a potentially asymmetric scaling of the viscosity near Dirichlet boundaries via
    \begin{flalign*}
		\hat{S}_{w}^{-1} \approx K_{1/\sqrt{\eta_l}}^{-1} \left ( B {M}_{\sqrt{\eta_l}}^{-1} A {M}_{\sqrt{\eta_r}}^{-1} B^T  \right ) K_{1/\sqrt{\eta_r}}^{-1},
	\end{flalign*}
    where 
    \begin{flalign*}
		\eta_l\rk{x} &:= \begin{cases}a_l^2 \eta\rk{x}  ,& x\in D_\eta\\\eta\rk{x},&\text{otherwise}\end{cases}
        \\
        \eta_r\rk{x} &:= \begin{cases}a_r^2 \eta \rk{x},& x \in D_\eta \\\eta\rk{x},&\text{otherwise,}\end{cases}
	\end{flalign*}
    for some $a_l, a_r > 0$. Here, $D_\eta=\bigcup \gk{K\in\mathcal{T}_{h},\partial K\cap\Gamma_{\Surface}\neq\emptyset}$. The scaled vector masses ${M}_{\sqrt{\eta_r}}$ and ${M}_{\sqrt{\eta_l}}$ are solved up to a relative tolerance of $\tol_{\VectorMass}>0$ using a CG solver preconditioned by a GMG approach.

    We propose a V-cycle BFBT variation of the diag(A)-BFBT approximation
	\begin{flalign} \label{eq:schurV}
		\hat{S}_{V}^{-1} :=	\left ( B \hat{A}_{C}^{-1} B^T \right )^{-1} \left ( B \hat{A}_{C}^{-1} A \hat{A}_{C}^{-1} B^T  \right ) \left ( B \hat{A}_{C}^{-1} B^T \right )^{-1},
	\end{flalign}
	replacing $\diag\rk{A}$ with $\hat{A}_{C}$, representing a comparatively computationally cheap $A$ block GMG V-cycle using only $m_V = 1$ pre- and post smoothing step with a Chebyshev smoother of degree $\deg_V = 1$. On the coarse grid, we again solve up to a relative tolerance of $\tol_{\coarse}>0$.
    
    In order to apply a V-cycle BFBT approximation, it is still necessary to solve an approximate Schur complement operator of the form $B \hat{A}_C^{-1} B^T$. Since a relatively large tolerance $\tol_{\VBFBT} \approx 0.1$ is sufficient for a solver suitable for geophysical applications with realistic viscosity contrasts, we use a CG solver preconditioned by $M_{1/\eta}$ which can be efficiently solved up to this comparatively large relative tolerance.
    \subsubsection{Saddle Point Solver Structure} \label{subsec:FGMRESSolver}
\begin{figure}[htbp]
        \includegraphics[width=\textwidth]{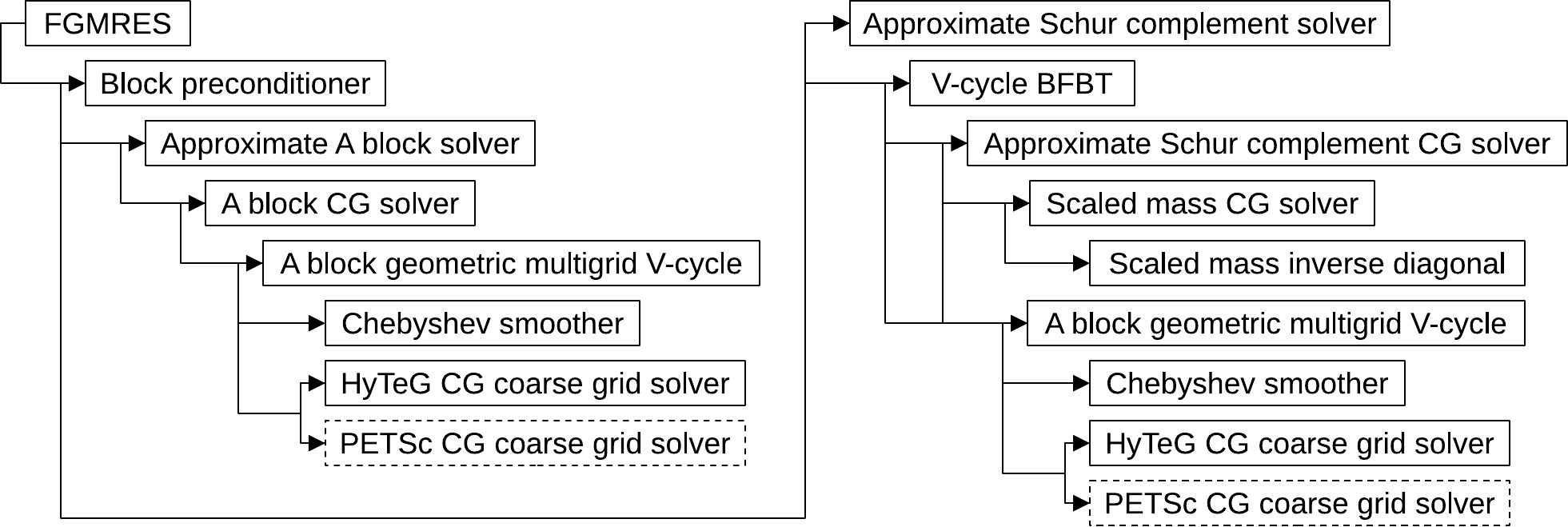}
    		\caption{Solver structure of the preconditioned FGMRES saddle point solver using the V-cycle BFBT Schur complement approximation.}
    		\label{fig:StokesSolverStructure}
    \end{figure}

	Altogether, our solver for the saddle point system \eqref{eq:SaddlePointStructure} consists of a flexible GMRES \cite{Saad1993} outer loop preconditioned with an Uzawa type block preconditioner as seen in Sec.~\ref{subsec:UzawaPreconditioner}. At each time step, we solve the saddle point system up to an absolute residual tolerance of $\tol_{\rk{u,p}}>0$. For small time steps, it is possible that this criterion is still fulfilled after updating the temperature. This allows for multiple temperature updates per saddle point solve (compare \cite[Remark 2.3]{Waluga2016}). Fig.~\ref{fig:StokesSolverStructure} shows the hierarchical structure of the solver in case of the V-cycle BFBT approximation.
    
	The combination of an outer Krylov solver loop and a block preconditioner is a well-known solver structure that has been successfully used in the past to solve Stokes problems (compare \cite{Kronbichler2012}, \cite{Heister2017}) and generalised Stokes problems \cite{Jodlbauer2024}, even in a high contrast setting \cite{Wichrowski2022}.

    \section{Discrete Temperature Problem and Solver} \label{sec:DiscTemp}
    Similar to the weak formulation of \eqref{eq:StokesSub}, dividing \eqref{eq:DiffSolve} by $\rho\rk{x} C^p$, multiplying by a test function $w \in H_0^1\rk{\Domain}$,
    integrating over the domain, using integration by parts, using the Gauß divergence theorem, and eliminating the boundary integral of the diffusion term due to the vanishing boundary of the test function, admits a variational problem of the form
    \begin{flalign}
        &\int_{\Domain} s^{n+1} T^{n+1} w\dx 
        
        + \tau^{n+1} \Bigg(
        
        \int_{\Domain} u_{*}^{n+1} \cdot \nabla T^{n+1} w\dx 

        + \int_{\Domain} \frac{k}{\Peclet C^p} \nabla T^{n+1} \cdot \nabla \rk{\frac{w}{\rho\rk{x}}} \dx

        \nonumber \\

        &- \int_{\Domain} \frac{\Dissipation \alpha}{C^p} T^{n+1} \rk{u_{*}^{n+1} \cdot g} w \dx

        \Bigg)

        = \tau^{n+1} \Bigg(

        \int_{\Domain} s^n \hat{T}^{n}  w_h \dx 

         - \int_{\Domain} s^{n-1} \hat{T}^{n-1}  w \dx 

         \nonumber \\

          & + \int_{\Domain} \frac{H}{C^p} w  \dx
        
         + \int_{\Domain} \frac{\Peclet \Dissipation}{\Rayleigh \rho\rk{x} C^p} 2 \eta\rk{x, T_{*}^{n+1}} \dot{\varepsilon}\rk{u_{*}^{n+1}}: \dot{\varepsilon}\rk{u_{*}^{n+1}}w  \dx

         \Bigg)
         \label{eq:WeakTemp}
    \end{flalign}
    \noindent
    for all $w \in H_0^1\rk{\Domain}$ with 
    \begin{flalign*}
         &s^{n+1} := \frac{2 \tau^{n+1} + \tau^{n}}{\tau^{n+1} + \tau^{n}}
         
         , \quad
         
         s^{n} := \frac{\tau^{n+1} + \tau^{n}}{\tau^{n}}
         
         , \quad
         
         s^{n-1} := \frac{\tau^{n+1} \cdot \tau^{n+1}}{\tau^{n} \left ( \tau^{n+1} + \tau^{n} \right )}.
    \end{flalign*}    
    \noindent
    Replacing $u_*^{n+1}\in V$, $T^{n+1}, T_*^{n+1}, \hat{T}^{n}, \hat{T}^{n-1} \in W$ in \eqref{eq:DiffSolve} by $u_{*,h}^{n+1}\in V_h$, $T_h^{n+1}$, $T_{*,h}^{n+1}$, $\hat{T}_h^{n}$, $\hat{T}_h^{n-1} \in W_h$ and the test function $w$ by 
    \begin{flalign*}
        w_h \in W_{0,h} := \Big\{ & g_h \mid g_h = \tilde{g}_h \circ B^{-1} \text{ for some } \tilde{g}_h \in C^0\rk{\tilde{\Domain}} \text{ with } 
        \\
        & \evalat{\tilde{g}_h}_{\rk{\partial \tilde{\Domain}}}  = 0, 
        \evalat{\tilde{g}_h}_{\tilde{K}} \in \mathcal{P}_2\rk{\tilde{K}} \, \forall \tilde{K} \in \tilde{\mathcal{T}}_h \Big\}
    \end{flalign*}
    leads to the discrete variational advection diffusion problem we solve in practice.
    
    Given $u_{*,h}^{n+1}\in V_h$ and $T_{*,h}^{n+1}, \hat{T}_h^{n}, \hat{T}_h^{n-1} \in W_h$, we seek to find an approximation $T_h^{n+1} \in W_h$, that fulfills the discrete variational advection diffusion problem for all test functions $w_h \in W_{0,h}$.
    
    To solve the linear system arising from the discrete variational advection diffusion problem, we use an FGMRES outer loop preconditioned by a CG solver solving the symmetric, positive semidefinite mass and diffusion part of the left-hand-side. The advection-diffusion system is solved up to an absolute tolerance of $\tol_{T}>0$. The application of the preconditioner is carried out up to machine precision in an exact way in each FGMRES iteration. The Dirichlet boundary conditions at the surface and CMB are again enforced in a similar manner to \eqref{eq:SaddlePointStructureDirichlet}.

	\section{Geodynamical Application and Numerical Results} \label{sec:GeodynamicApplication}
    In this chapter, we present several numerical tests showcasing the convergence rate in time, the scalability of the code, a comparison of convergence of all proposed solvers, and the applicability to real-world problems using high-resolution simulations.
    \subsection{Temporal Discretisation Convergence Test}
    We test the numerical convergence of the time stepping scheme by using \eqref{eq:DiffSolve} with $\rho\equiv1=C^p=\Dissipation=\alpha=H=\Peclet=\Rayleigh$, resulting in
    \begin{equation}\label{eq:convergenceTest}
\frac{\partial T}{\partial t} + \left(u \cdot \nabla\right) T - k \Delta T -  T  \left ( u \cdot g \right ) -  \tau : \dot{\varepsilon} - 1 =  f,
    \end{equation}
    where $k\in\{10^{-5},10^{-4},\ldots,10^0,3\}$ and $f$ will be explained later. The test is run in two-dimensional space with $r_{\CMB}=\frac{1}{2}$ and $r_{\Surface}=\frac{3}{2}$. As time interval we use $[\frac{7}{2},\frac{9}{2}]$. We set 
    \[
    \begin{split}
        T(t,x) &= x_0 x_1 \cos\rk{ t \sqrt{\rk{x_0+\sin(3t)}^2 + \rk{x_1-\cos(3t)}^2} },\\
        u(t,x)&=(2+\cos(3\pi t+x_0 x_1 t))\begin{pmatrix}-x_1\\x_0\end{pmatrix},\text{ and}\\
        \eta(x,T)&=\norm{x}^2 \exp\rk{-T}
    \end{split}
    \]
    as solution for the temperature, the velocity, and as temperature-dependent viscosity. Given these solutions, we solve \eqref{eq:convergenceTest} according to the schemes \eqref{eq:timeDerivative}--\eqref{eq:Extrapolation}, where the previously mentioned right-hand side $f$ is the left-hand side of \eqref{eq:convergenceTest} evaluated with the solutions at a given time point. We use an equidistant time grid. For the first time step of the multi-step method, we evaluate the solutions for the needed time steps before the initial time point.
    \begin{figure}[htbp]
		\centering
		\includegraphics[height=0.22\textheight]{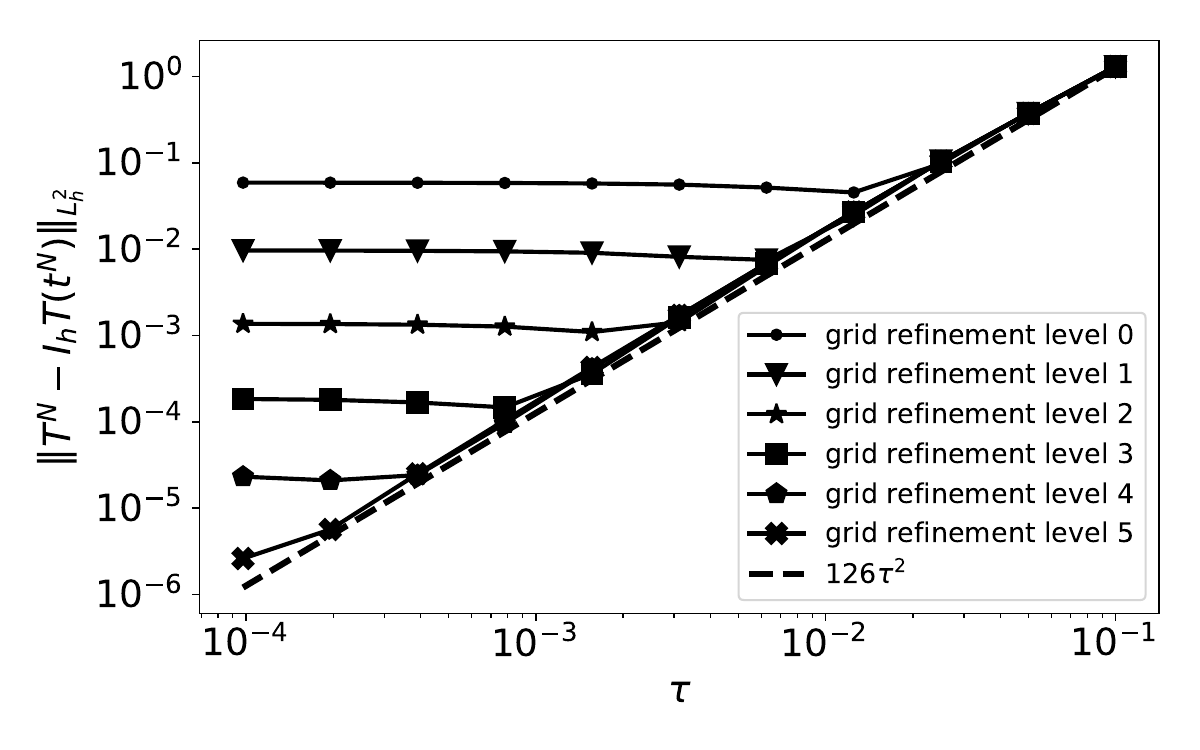}
        \includegraphics[height=0.22\textheight]{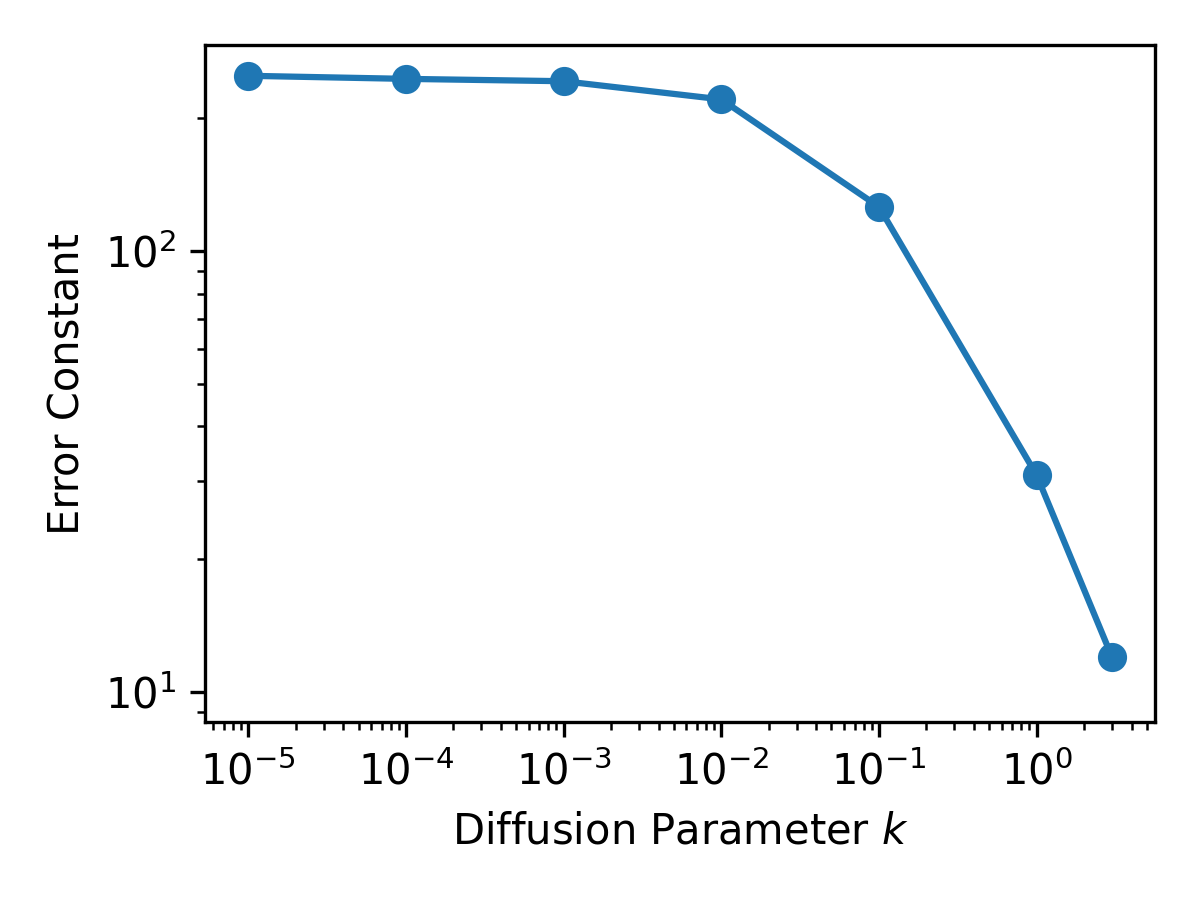}
		\caption{Absolute $L^2$-error of the approximated numerical solution at time point $t^{N}$ for different grid refinement levels of the blended unstructured annulus coarse grid plotted against the time step size $\tau$ for $k=10^{-1}$ (left). Error constant (in front of $\tau^2$) depending on the diffusion parameter $k$ (right).
        } \label{fig:temporalDiscResults}
	\end{figure} 	    

    The left of Fig.~\ref{fig:temporalDiscResults} depicts the discrete $L^2$-error of the final time step, i.e., $\left \|T^{N} - I_hT\left (\frac{9}{2}\right ) \right \|_{L^2_h}$, where $I_h$ is a projection onto the spatial grid, for different coarse grid refinement levels plotted against $\tau=N^{-1}$. The time discretisation scheme leads to a second order temporal convergence given a fine enough spatial discretisation. The diffusion parameter $k$ in \eqref{eq:convergenceTest} does not impact the order of convergence, but the constant in front of $\tau^2$. The right of Fig.~\ref{fig:temporalDiscResults} shows that while the constant increases with decreasing $k$, the constant seems to be bounded. This suggests, that both, MMOC and the diffusion solve, have a second order convergence rate in time, with the MMOC handling advection-dominated problems having a larger but bounded error constant.
	\subsection{Time Dependent Geodynamical Model} \label{subsec:GeoModel}
    As described in Sec.~\ref{sec:Model}, we assume $\alpha$, $k$, $C^p$ and $H$ to be constant, $g(x) = -\frac{x}{\norm{x}_2}$, the density $\rho\rk{x}$ to be space dependent and the viscosity $\eta\rk{x,T}$ to be space and temperature dependent. Tab.~\ref{dimtable} shows our choice for $\Rayleigh$, $\Peclet$ and $\Dissipation$ along the chosen reference constants for nondimensionalisation. In particular, we assume $\alpha$, $k$, $C^P$ and $H$ to be equal to their reference constants, hence we choose $\alpha=k=C^P=H=1$ in the nondimensional equations.

\begin{figure}[htbp]
		\centering
		\includegraphics[width=0.7\textwidth]{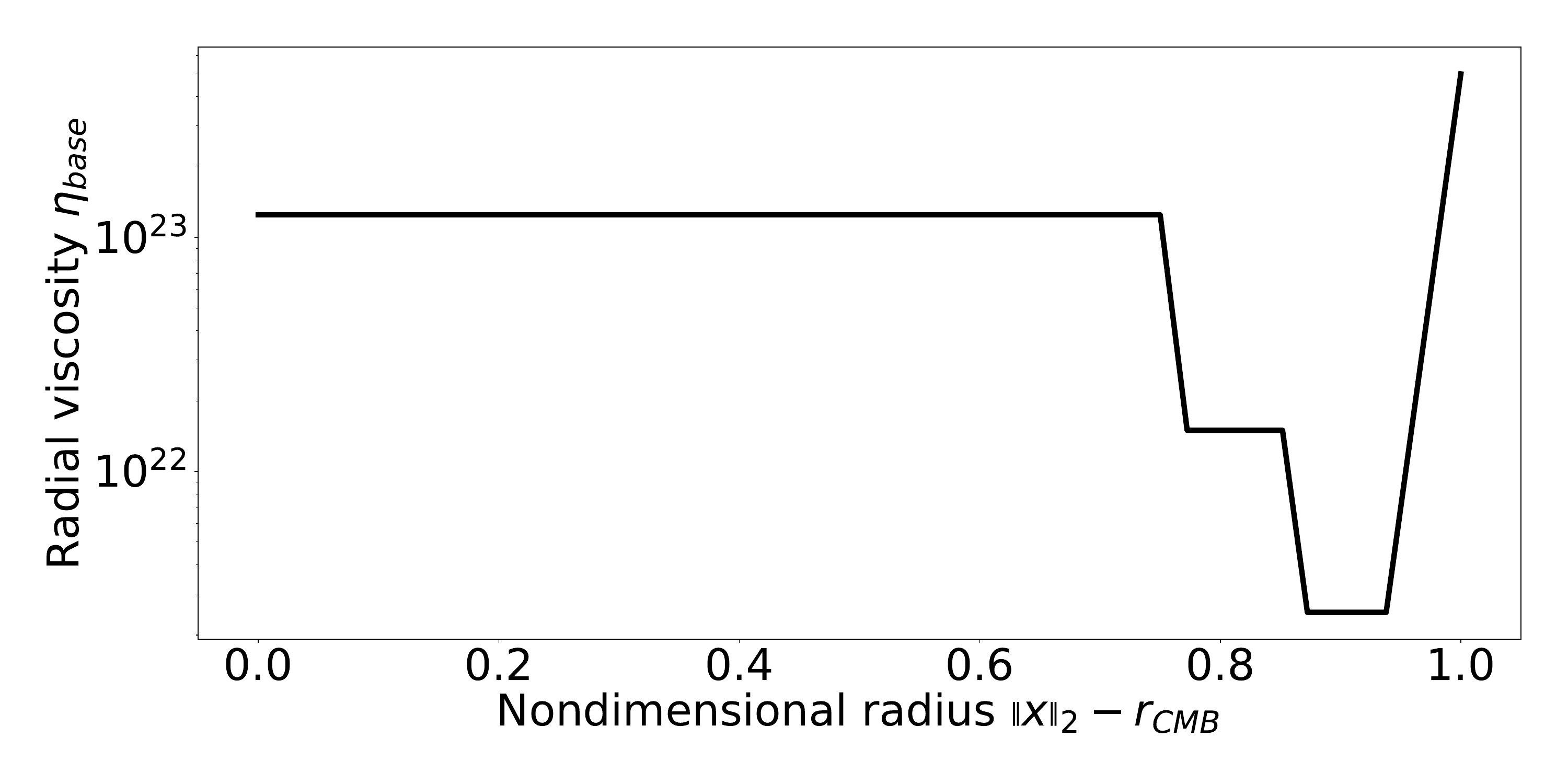}
        \caption{Radial viscosity $\eta_{\base}$ plotted against the nondimensional radius $\norm{x}_2 - r_{\CMB}$. This profile is similar to the viscosity contrast shown in \cite[Fig.~8]{Lin2022}.}
        \label{fig:baseViscosity}
	\end{figure}

    We choose a temperature dependent Frank–Kamenetskii type viscosity \cite{May2015, Lin2022}
	\begin{equation}
		\eta(x,T) = \eta_{\base} \left ( x \right ) \exp\left ( -E_A \cdot T + V_A \cdot \left ( r_{\Surface} - \left \| x \right \|_2 \right ) \right ) \label{eq:ViscosityDef}
	\end{equation}
     with activation energy $E_A := 4.610$ and activation volume $V_A := 2.996$ (compare \cite{Lin2022}) modulating the radial base viscosity profile $\eta_{\base}$ depicted in Fig.~\ref{fig:baseViscosity}.   
     \noindent
	For the density and reference temperature, we choose
	\begin{align*}
		\rho\left ( x \right ) = \rho_{\Top}  \exp\left ( \frac{\Dissipation}{\Gamma_0} \left ( r_{\Surface} - \left \| x \right \| \right )  \right )
	      \\ 
		T^s = T_{\text{adiabatic}} \exp\left ( \Dissipation \left ( r_{\Surface} - \left \| x \right \| \right )  \right )
	\end{align*}    
	with $\rho_{\Top} = 3.381 \cdot 10^{3} \ukg \um^{-3}$ and $T_{\text{adiabatic}} = 1.6 \cdot 10^{3} \uK$. As an initial temperature, we use $T^s$ plus a 3\% relative uniform random noise.

    As a Schur complement approximation, we choose the V-cycle BFBT method. As a block preconditioner for the saddle point system, we choose the symmetric Uzawa block preconditioner. Tab.~\ref{parametertable} lists the default solver tolerances and parameters chosen in the following numerical simulations unless stated otherwise.

    To avoid unphysical interactions with the plate velocities prescribed at the surface, the shear heating term has been set to zero close to the surface for the following tests.
\begin{figure}[htbp]
		\centering
        \includegraphics[width=0.5\textwidth]{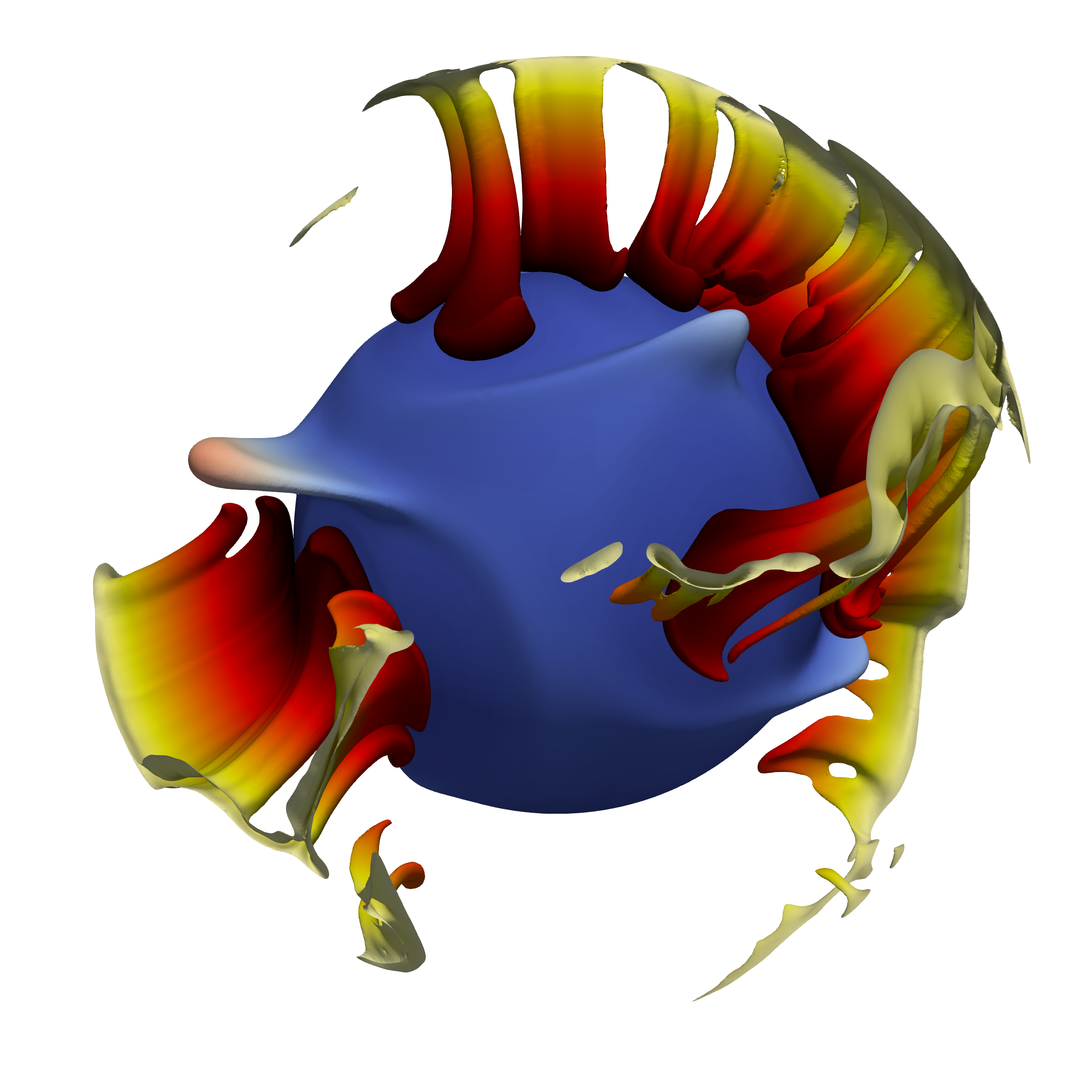}\includegraphics[width=0.5\textwidth]{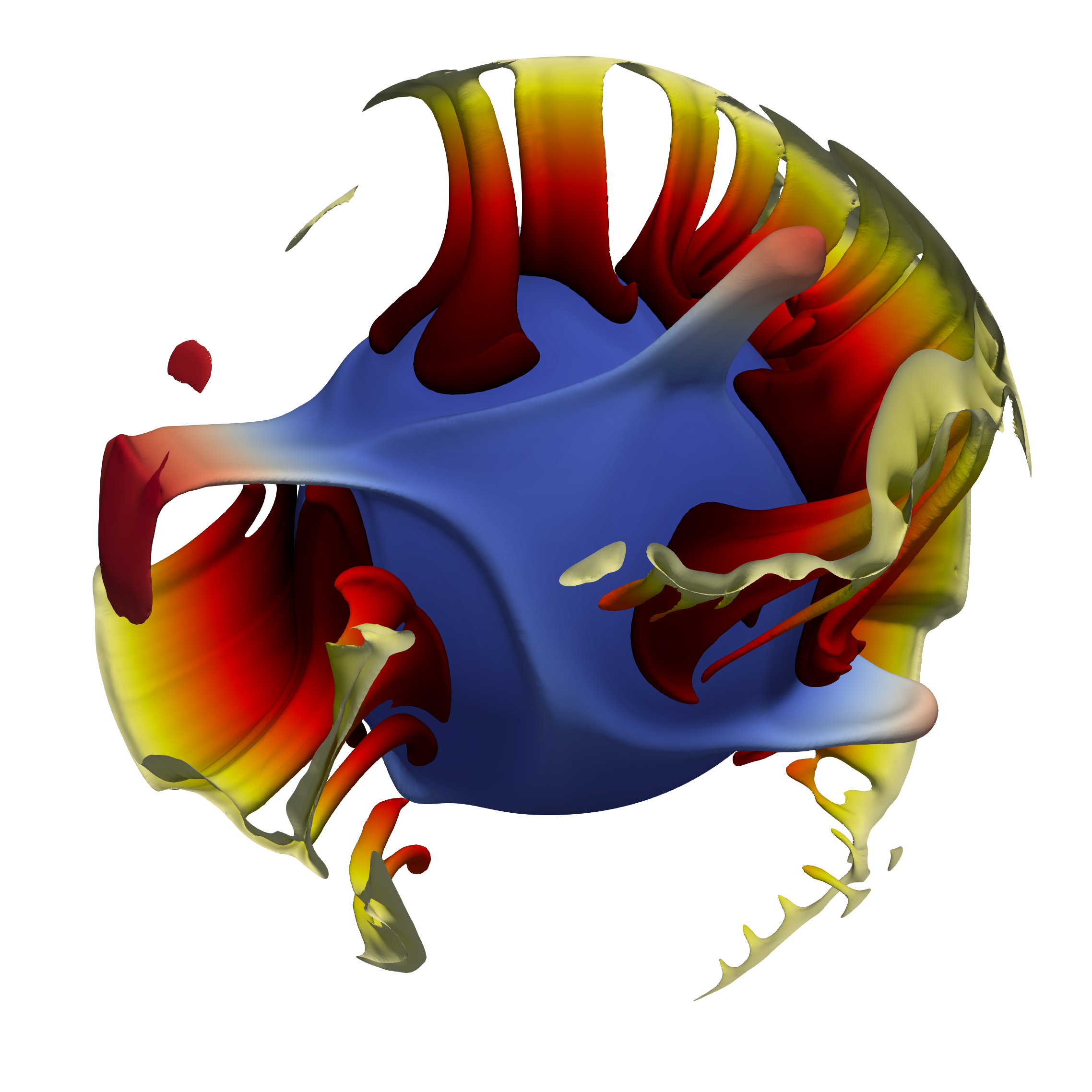}
        \\
        \includegraphics[width=0.49\textwidth]{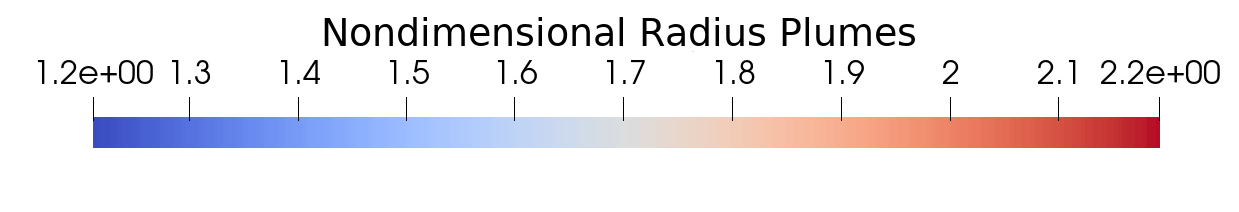}
        \hfill
        \includegraphics[width=0.49\textwidth]{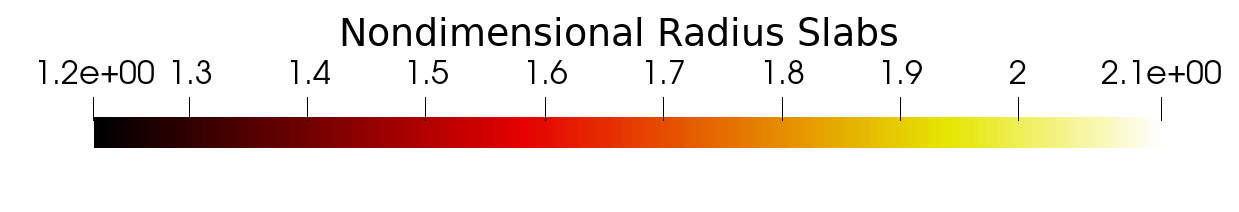}
        \caption{Simulation results at 175 Myr (left) and 200 Myr (right), depicting contour plots of the temperature deviation, linearly colored by the distance to the CMB for plumes (blue to red) and slabs (dark red to yellow). One can clearly see how three plumes form, one reaching the surface of the Earth (mid left in both pictures).}
        \label{fig:175Myr_Base}
    \end{figure}   

    \begin{figure}[htbp]
		\centering
        \includegraphics[width=0.85\textwidth]{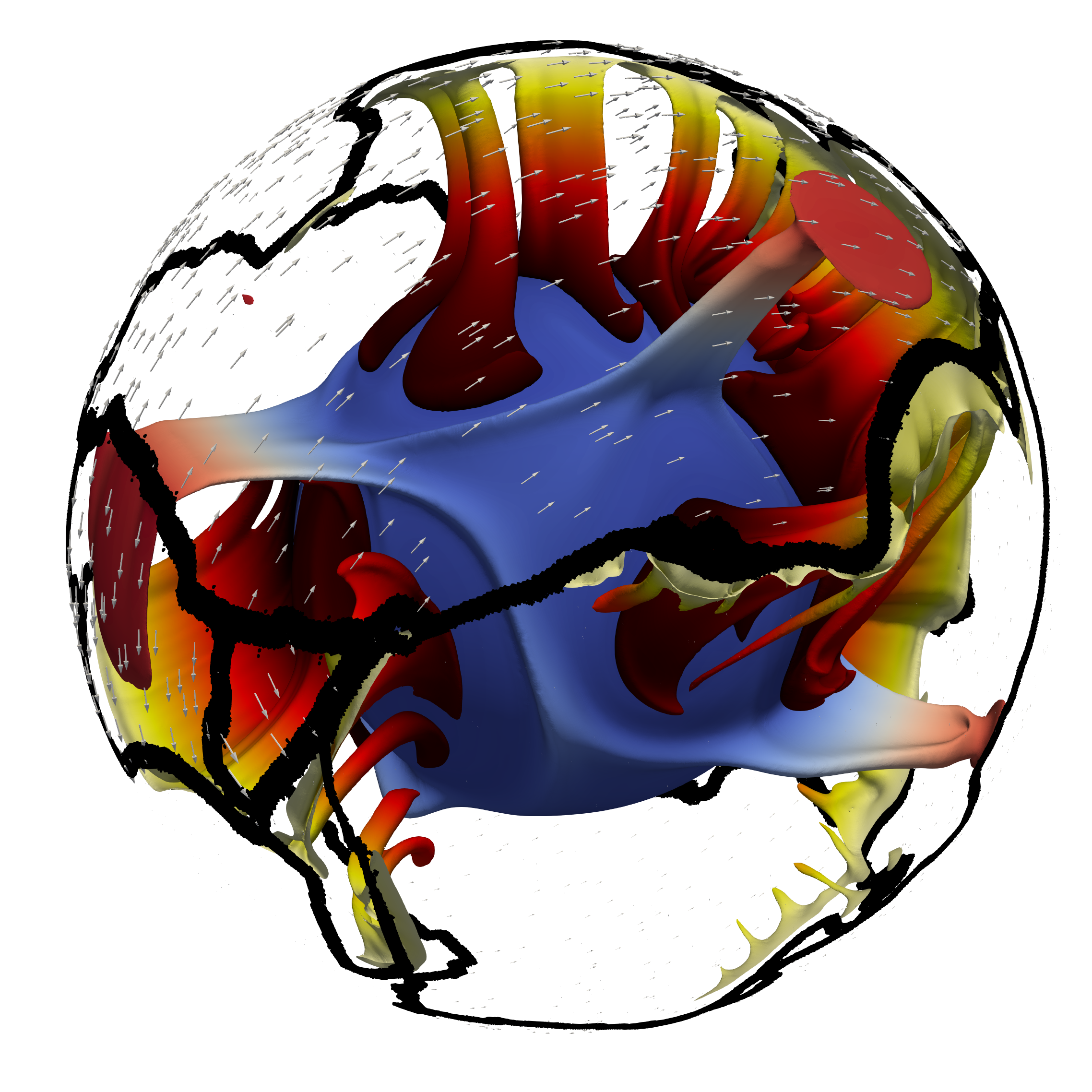}
        \caption{Simulation result at 212.5 Myr, depicting contour plots of the temperature deviation, linearly colored by the distance to the CMB for plumes (blue to red) and slabs (dark red to yellow). The black area on the surface shows the boundary of the tectonic plates. Grey arrows indicate the movement of each plate. While one plume (mid left) is displaced heavily by the plate movement, another plate (top right) is displaced much less.}
        \label{fig:212Myr_Plates_Vel}
    \end{figure}     

    \begin{figure}[htbp]
		\centering
        \includegraphics[height=0.55\textwidth]{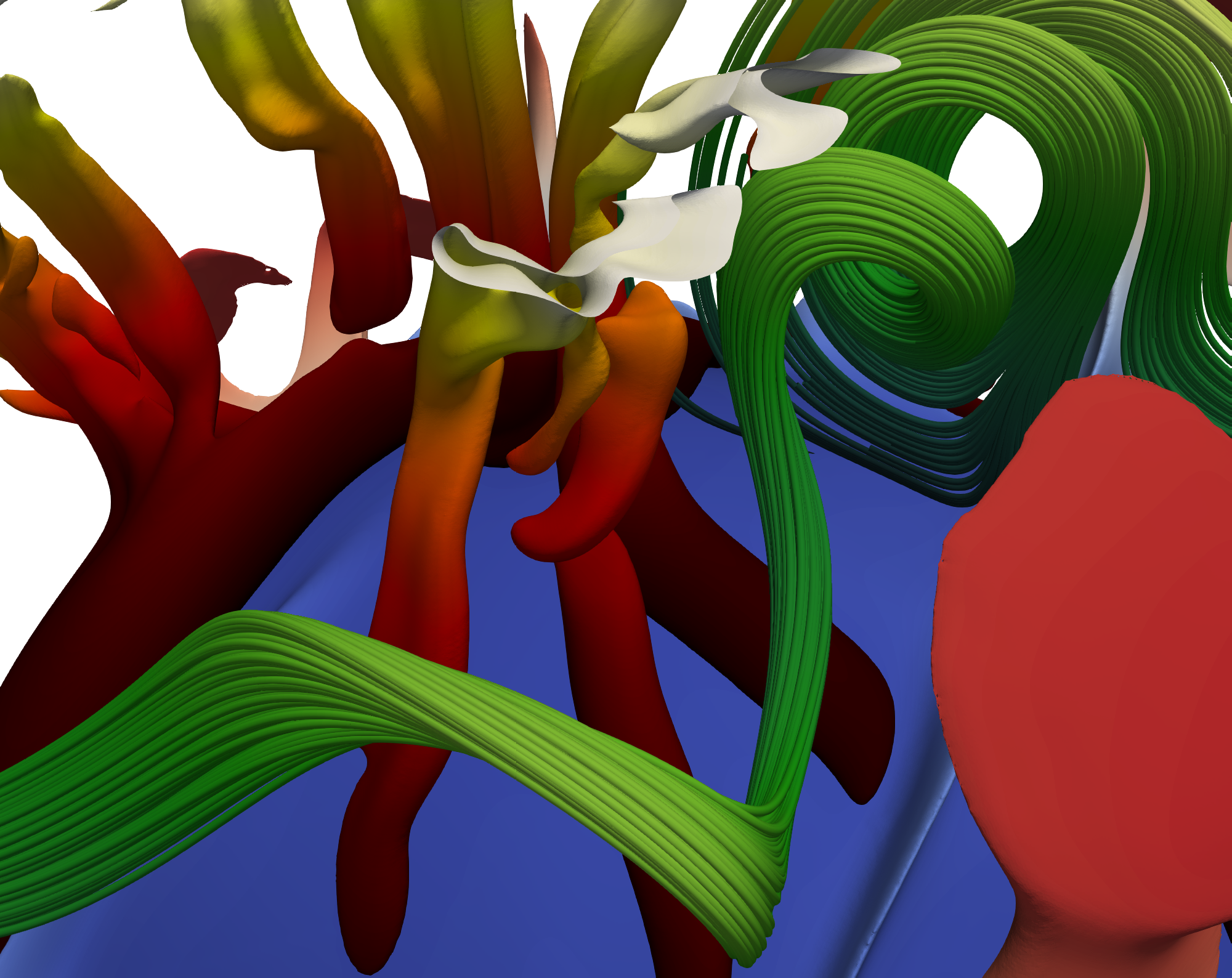}
        \includegraphics[height=0.55\textwidth]{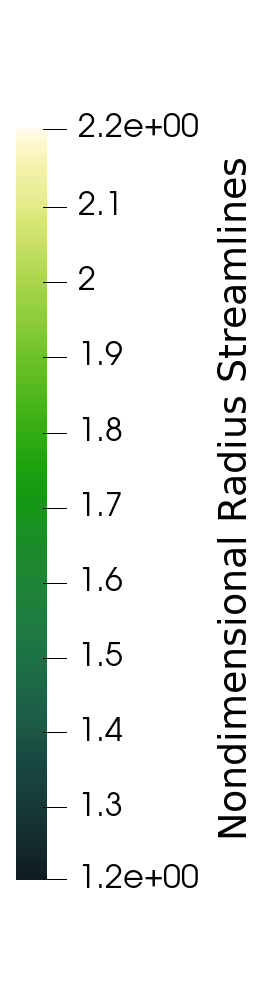}
        
        \caption{Simulation result at 312.5 Myrs, showcasing streamlines leading up to a convection cell within the Earth's mantle, linearly colored by their distance to the CMB (green to light yellow). The same coloring is used as in the previous figures for slabs and plumes.}
        \label{fig:Streamlines}
    \end{figure}
\begin{figure}[htbp]
		\centering
        \includegraphics[width=0.328\textwidth]{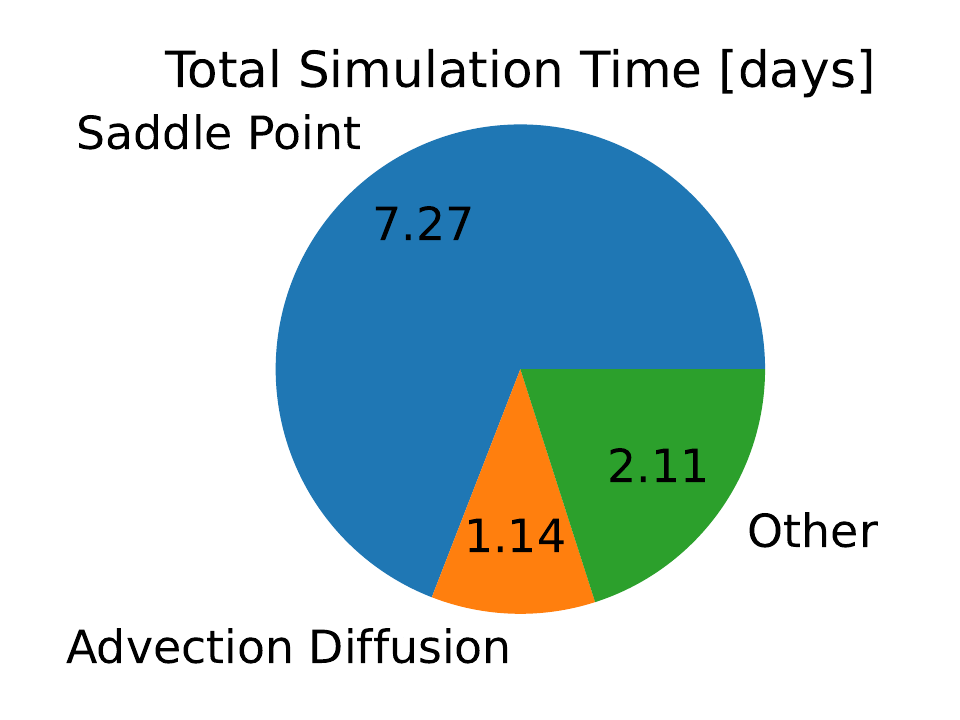}
        \hfill
        \includegraphics[width=0.323\textwidth]{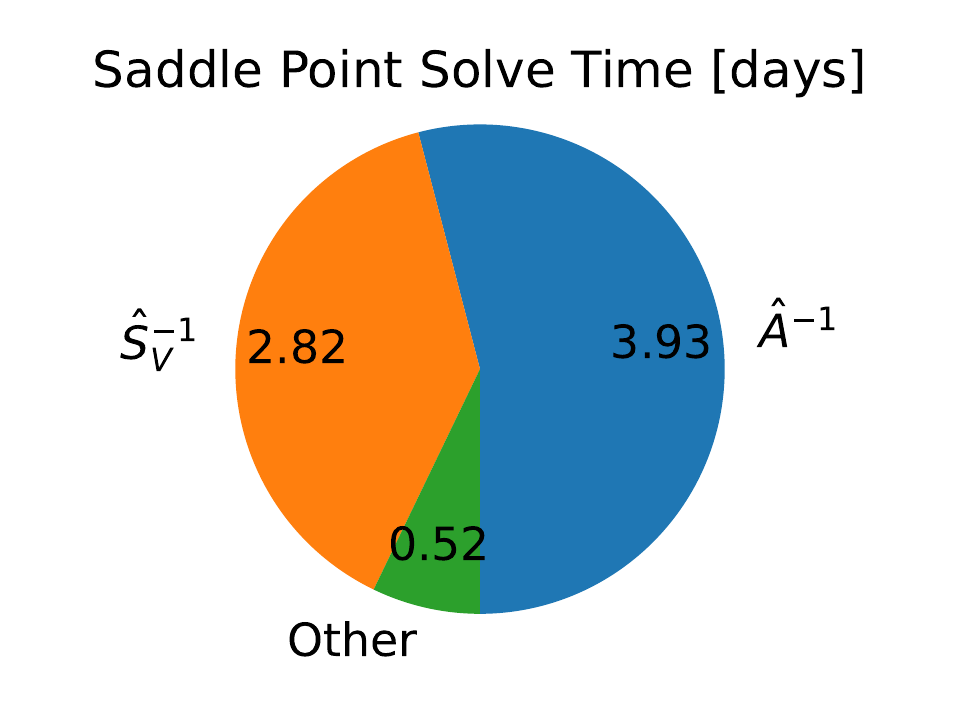}
        \hfill
        \includegraphics[width=0.328\textwidth]{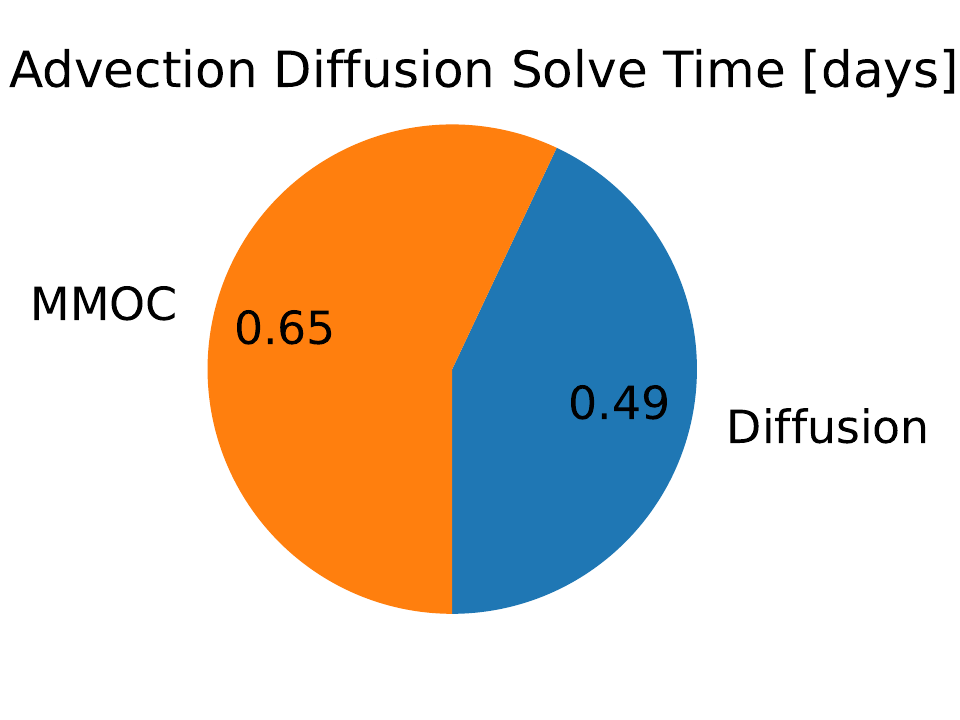}
        \caption{Time spent for the overall simulation (left), for solving the saddle point problem (middle) and the advection diffusion problem (right). The \textit{other} part for the total simulation time also includes calculating the time steps, reading the plate data, pre-calculating the viscosity interpolation, writing the output and calculating right-hand sides and extrapolations. The \textit{other} part for the saddle point problem includes the static part of the FMGRES and its block preconditioner and handling helper functions.}\label{fig:PIECHART}
    \end{figure}
\begin{figure}[htbp]
		\centering
        \includegraphics[width=0.5\textwidth]{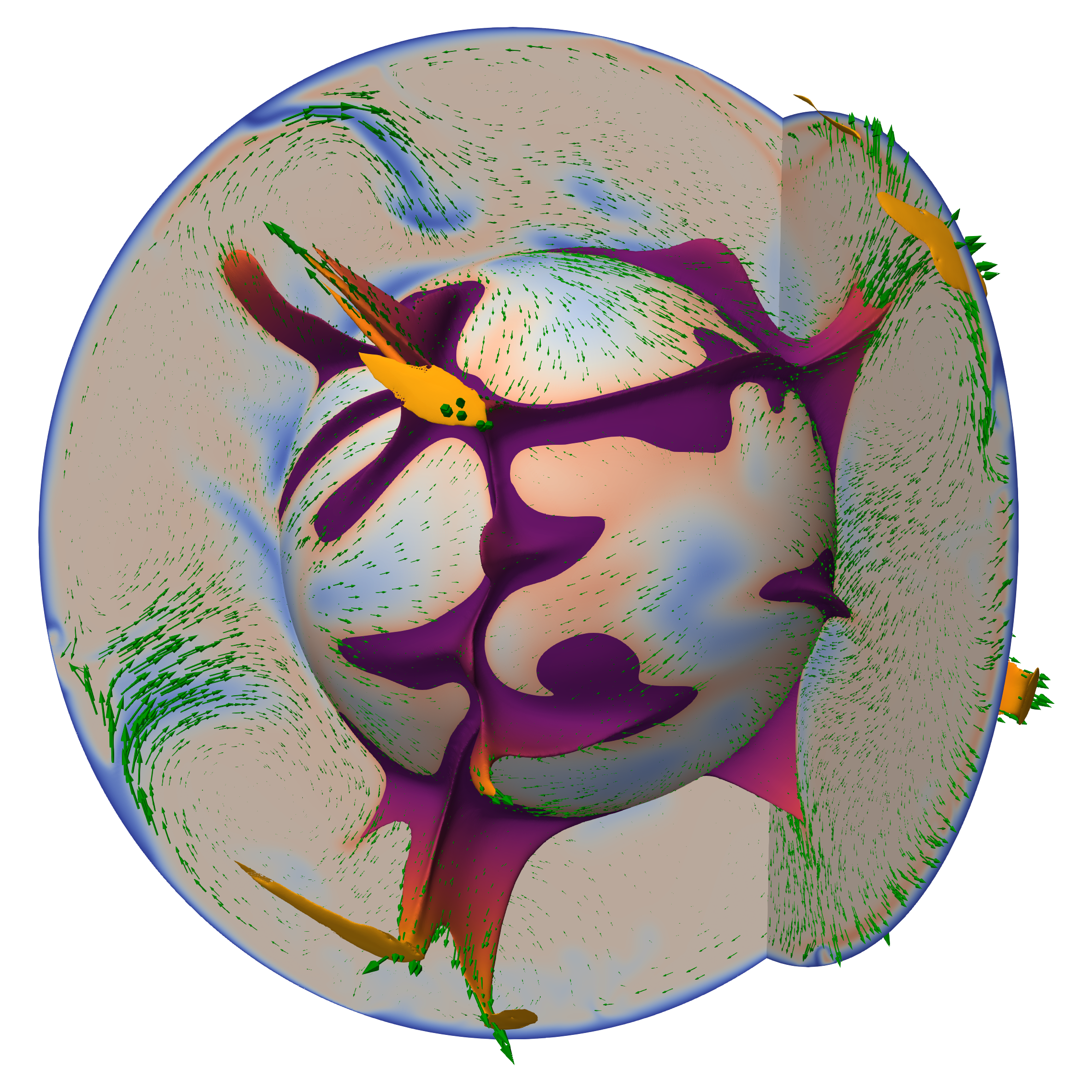}\includegraphics[width=0.5\textwidth]{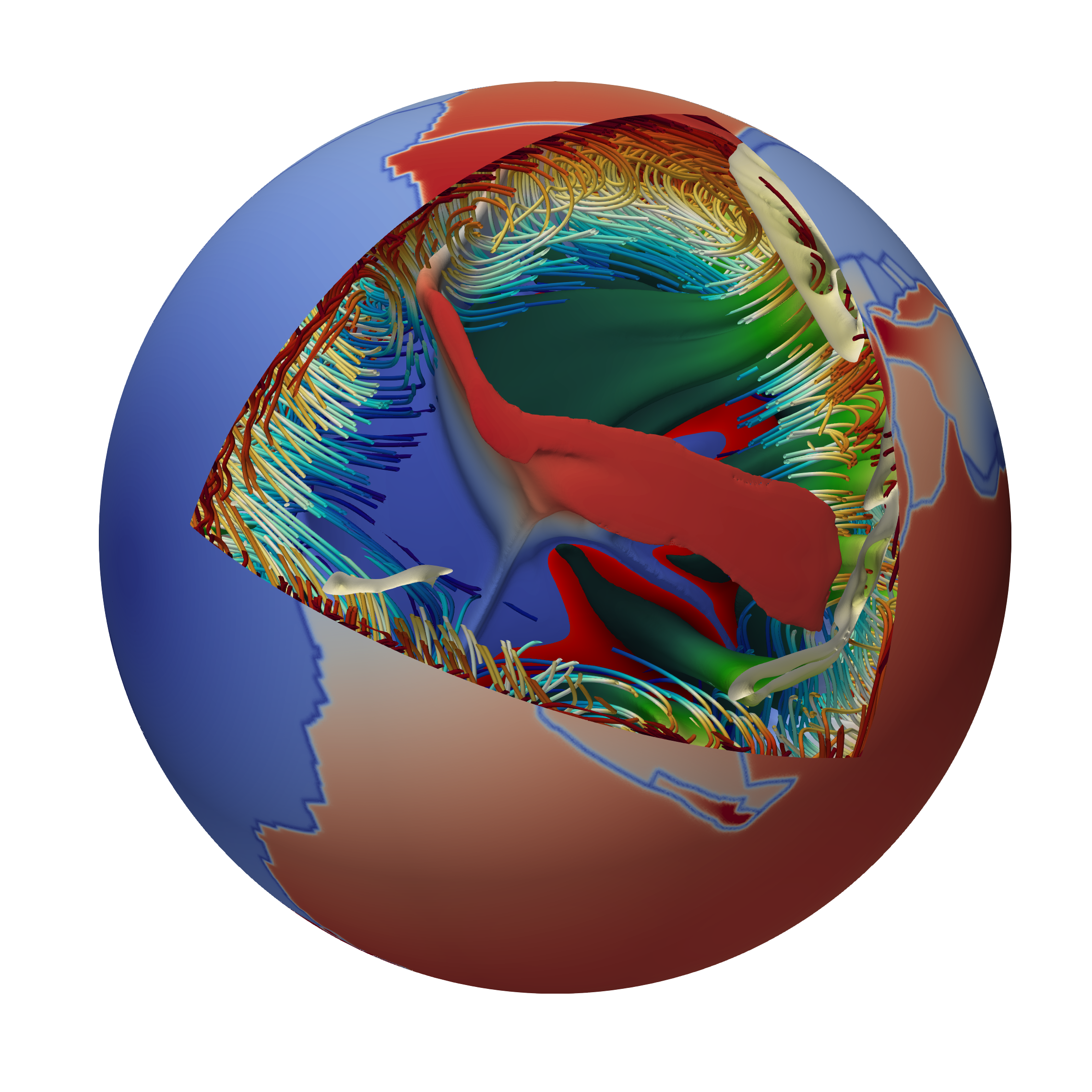}
        \\
        \includegraphics[width=0.49\textwidth]{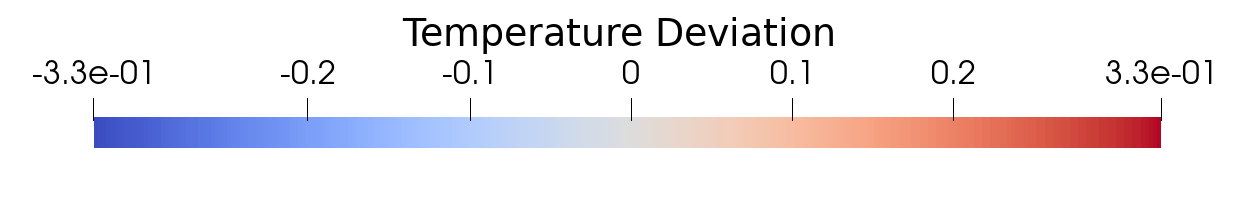}
        \hfill
        \includegraphics[width=0.49\textwidth]{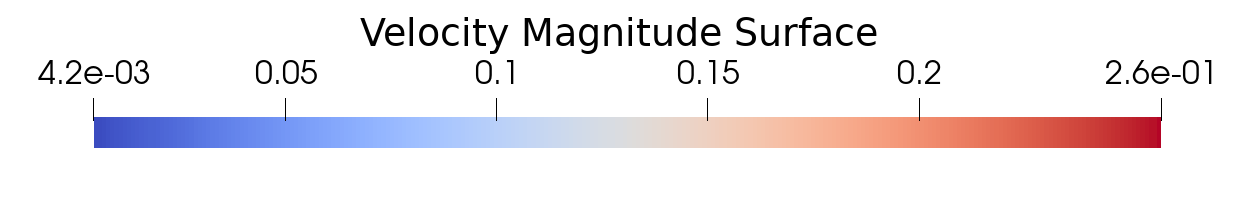}
        \\
        \includegraphics[width=0.49\textwidth]{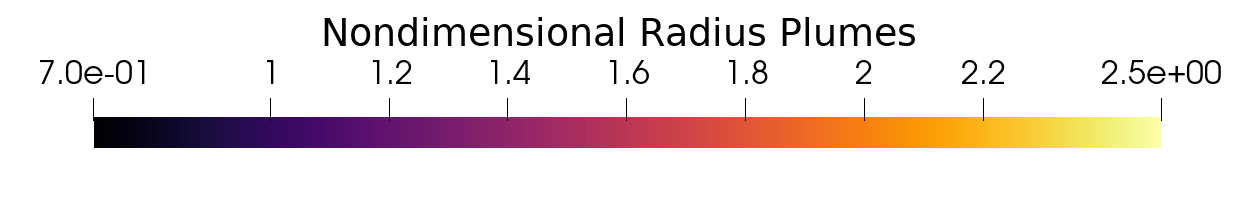}
        \hfill
        \includegraphics[width=0.49\textwidth]{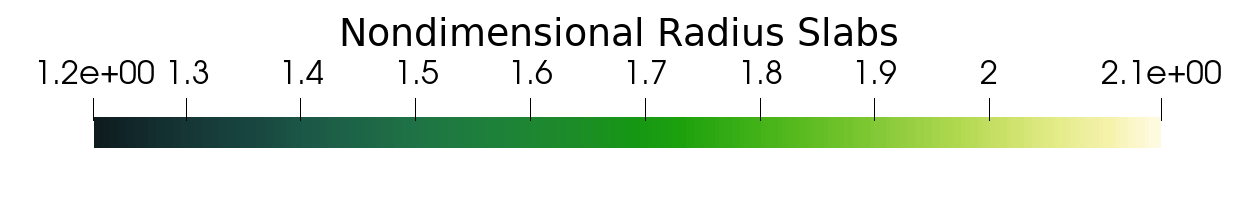}
        \\       
        \hfill
        \includegraphics[width=0.49\textwidth]{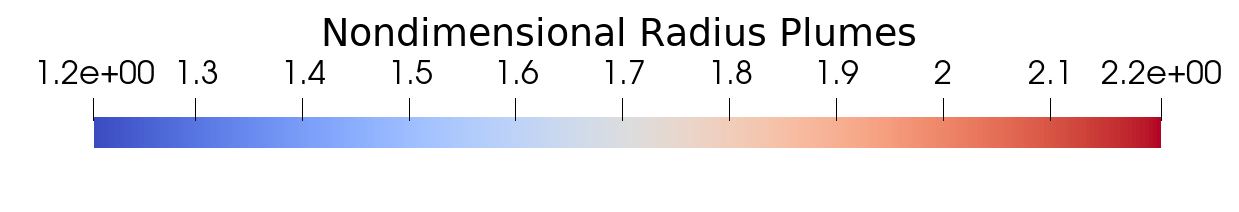}
        \\          
        \hfill
        \includegraphics[width=0.49\textwidth]{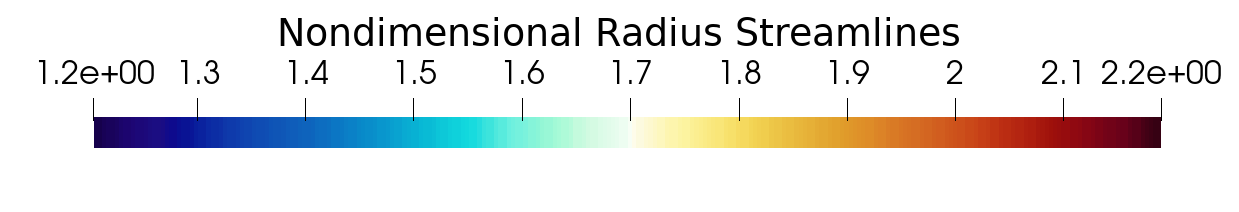}
        \\      
        \hfill
        \includegraphics[width=0.49\textwidth]{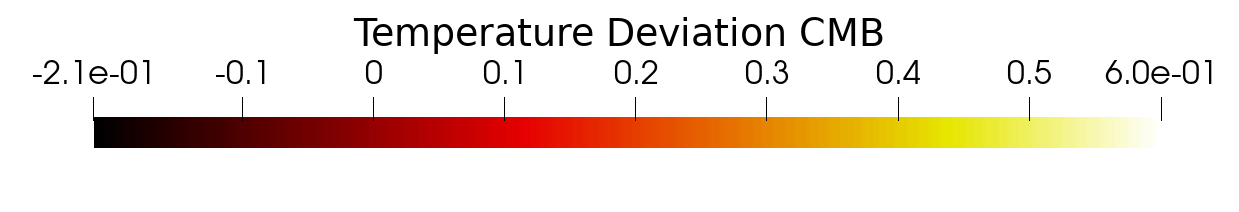}
        \\            
        \begin{minipage}{0.49\textwidth}
        \begin{spacing}{0.6}
        {\footnotesize \textbf{(a)} On the CMB and two slices, the temperature deviation is depicted from cold (blue) to warm (red) while the velocity field is shown via green arrows. Furthermore, plumes are visualized via a contour plot of the temperature deviation colored purple to orange, depending on the distance to the Earth's core.}
        \end{spacing}
        \end{minipage}
        \hfill
        \begin{minipage}{0.49\textwidth}
        \begin{spacing}{0.6}
        {\footnotesize \textbf{(b)} On the surface of the Earth, the plate data is shown from low velocity (blue) to high velocity (red). Slabs are visualized in green (dark to light, depending on the distance to the Earth's core), plumes in blue to red. Stream lines indicate the velocity field, linearly colored from blue to yellow and red, also depending on the distance to the Earth's core. While convection cells are clearly visible in the upper part of the picture, the stream lines follow the movement of slabs and plumes.}
        \end{spacing}        
        \end{minipage}   
        \caption{Simulation result at 345 Myrs (left) and 257.5 Myrs (right).}
        \label{fig:2DSlices_3DPlumes}
    \end{figure}

    Using a high resolution in space, we have $2.106\cdot10^9$ degrees of freedom in $(u,p)$ and $6.737\cdot10^8$ in $T$. Going forward in time, we start with reconstructed plate velocities on the surface at 1000 Myrs in the past. During the simulation the plate velocities are updated every 1 Myrs. The simulation time is $400$ Myr (plate time) with 6149 time steps, the computation time was $10.52$ days on the supercomputer HAWK at the HLRS Stuttgart, using 15360 cores on 128 nodes, where each node has two sockets consisting each of an AMD EPYC 7702 64-Core Processor with 2.25 GHz, and 128GB of memory per socket. The computation time can be divided into $7.27$ days for solving the saddle point problem via FMGRES, of which $3.93$ days are for solving the $A$ block via CG and $2.82$ days are for the Schur complement solver, and into 1.14 days for the advection-diffusion problem ($0.49$ days for the diffusion and $0.65$ days for the MMOC), see Fig.~\ref{fig:PIECHART} for a summary of the computation time. Different scenarios from this run are shown in Figures \ref{fig:175Myr_Base}, \ref{fig:212Myr_Plates_Vel}, \ref{fig:Streamlines} and \ref{fig:2DSlices_3DPlumes}.

    \subsection{$A$ Block Solver Test} \label{sec:ABlockTest}
    In this section, we compare different choices for $l_{\eta}$ and $l_{\min}$ with respect to $\hat{A}$ as discussed in Sec.~\ref{subsec:AApproximation} in the setting described in Sec.~\ref{subsec:GeoModel} with $2.64 \cdot 10^8$ of degrees of freedom in $(u,p)$ and $8.45 \cdot 10^7$ degrees of freedom in $T$ on a machine with 240 cores (2 Sockets with AMD EPYC 9754 128-Core Processor @ 2.25 GHz each) for a coarse grid with 240 macro elements and a maximum refinement level $l_{\max}=6$. For this purpose, we observe the total number of V-cycle preconditioned $A$ block outer CG iterations and the total time spent solving the $A$ block for 10 FGMRES iterations applied to the initial Stokes solve.
    \begin{figure}[htbp]
		\centering
        \includegraphics[width=0.49\textwidth]{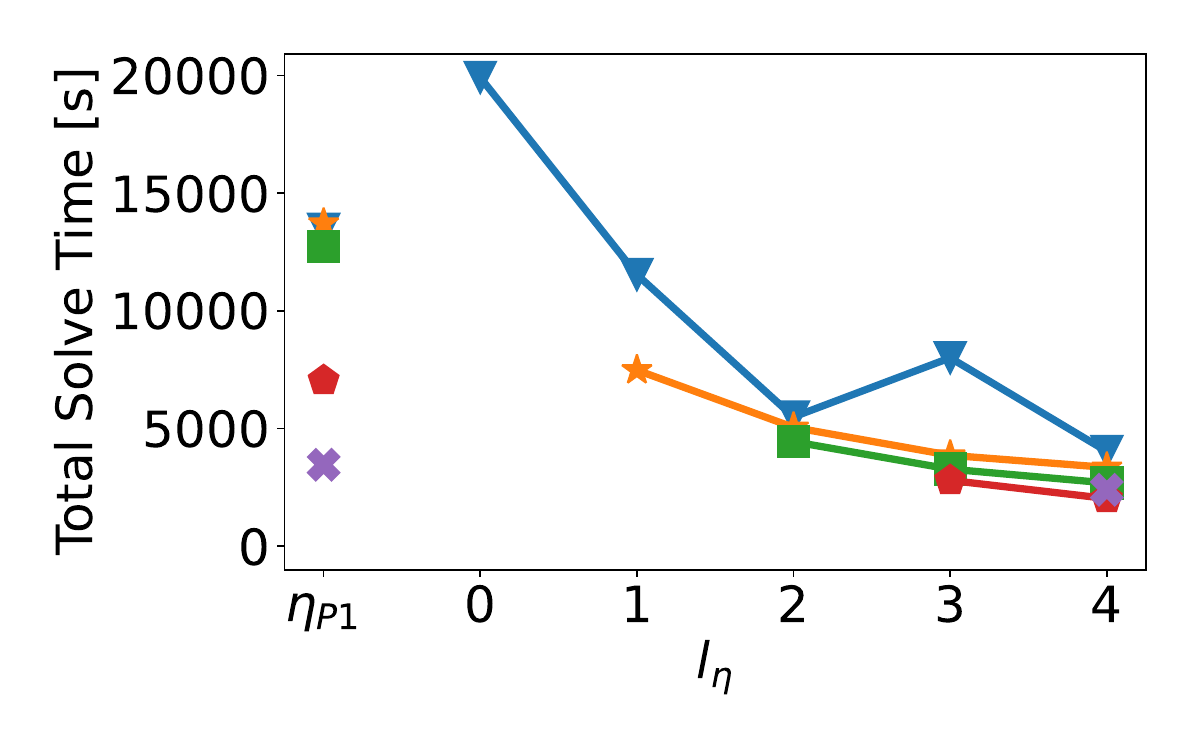}
        \includegraphics[width=0.49\textwidth]{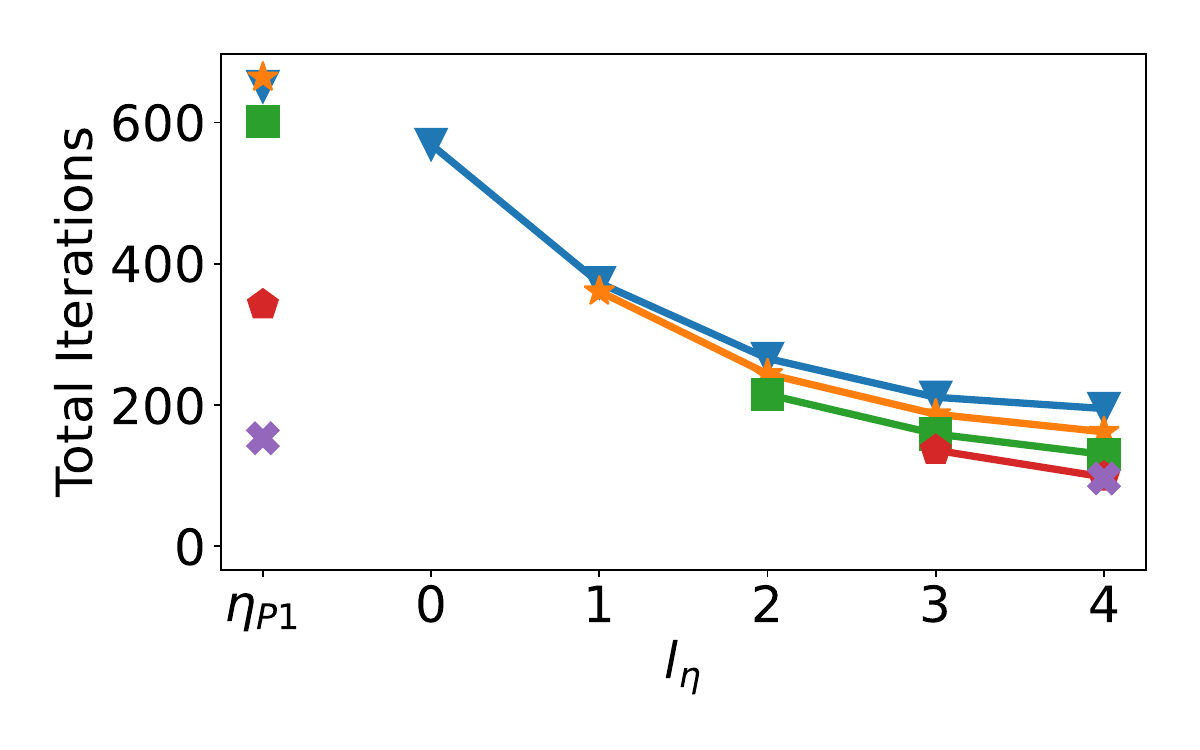}\\
        \includegraphics[width=0.49\textwidth]{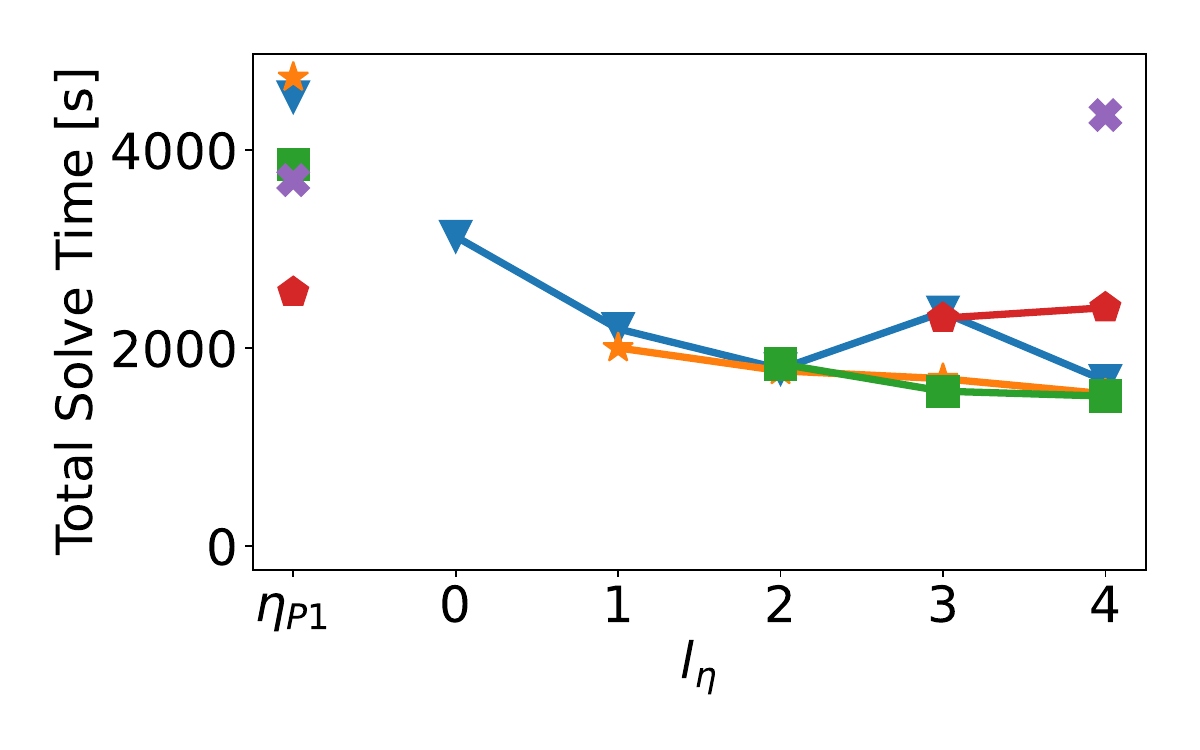}
        \includegraphics[width=0.49\textwidth]{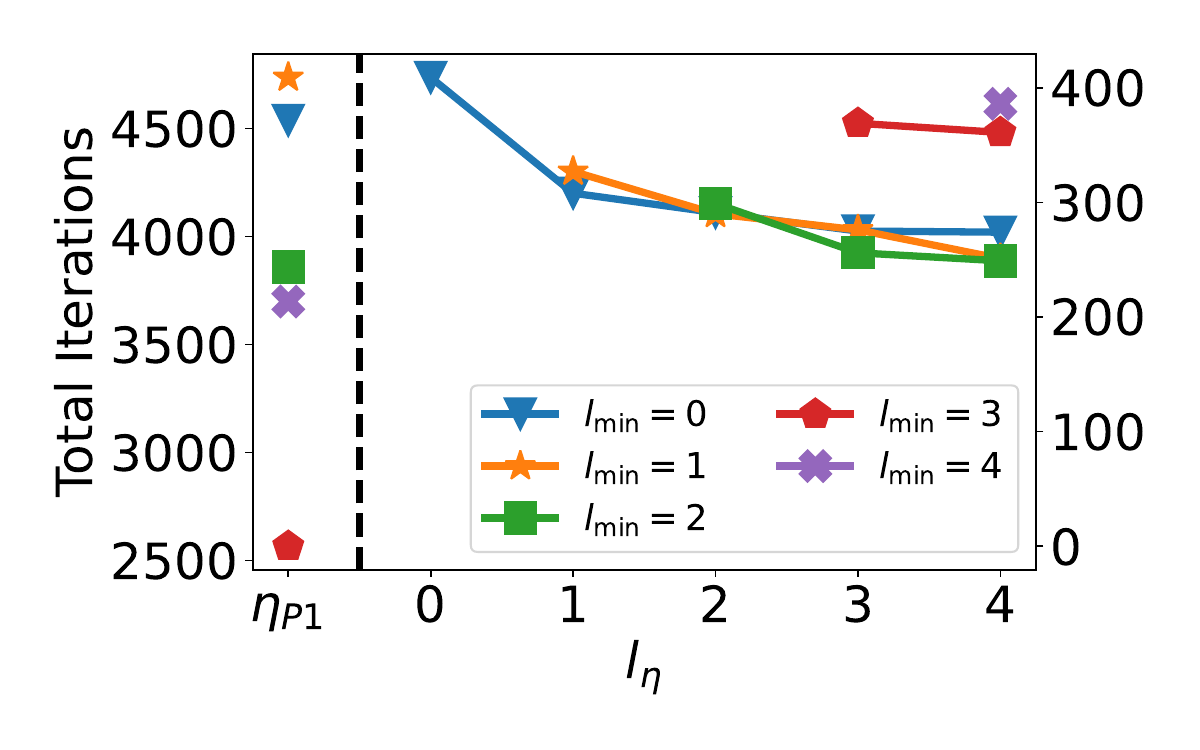}
        \caption{Total time spent (left) and total CG iterations (right) required to solve the $A$ block up to $\tol_A$ (top) or $B \hat{A}_C^{-1}B^T$ up to $\tol_{\VBFBT}$ (bottom) in 10 FGMRES iterations applied to the initial Stokes solve. The values at $\eta_{\text{P1}}$ denote evaluating $\eta^{n+1}_{\text{P1}}$ at the quadrature points on every refinement level. Note the split in axes for the bottom right figure.}
        \label{fig:AIterations}
    \end{figure} 
    
    Fig.~\ref{fig:AIterations} show the total number of iterations and total time spent solving the $A$ block up to $\tol_A$ as well as the Schur complement approximation operator $B \hat{A}_C^{-1}B^T$ up to $\tol_{\VBFBT}$, respectively. In practise, we choose $l_{\min}=2, l_\eta=3$ on a coarse grid with 240 macro elements, $l_{\min}=1, l_\eta=2$ on a coarse grid with 1920 macro elements, and $l_{\min}=0, l_\eta=1$ on a coarse grid with 15360 macro elements for our experiments. For general $\eta$, the choice of $l_\eta$ is delicate, as it has to balance between the total solve time and the memory needed for an implementation using matrices which sizes depend on the coarse grid. As in our case, for a specific $\eta$, when using an optimized matrix free implementation, larger $l_\eta$ are generally more desirable.

    \subsection{Schur Complement Approximation Comparison}
    \begin{figure}[htbp]
		\centering
        \includegraphics[width=0.49\textwidth]{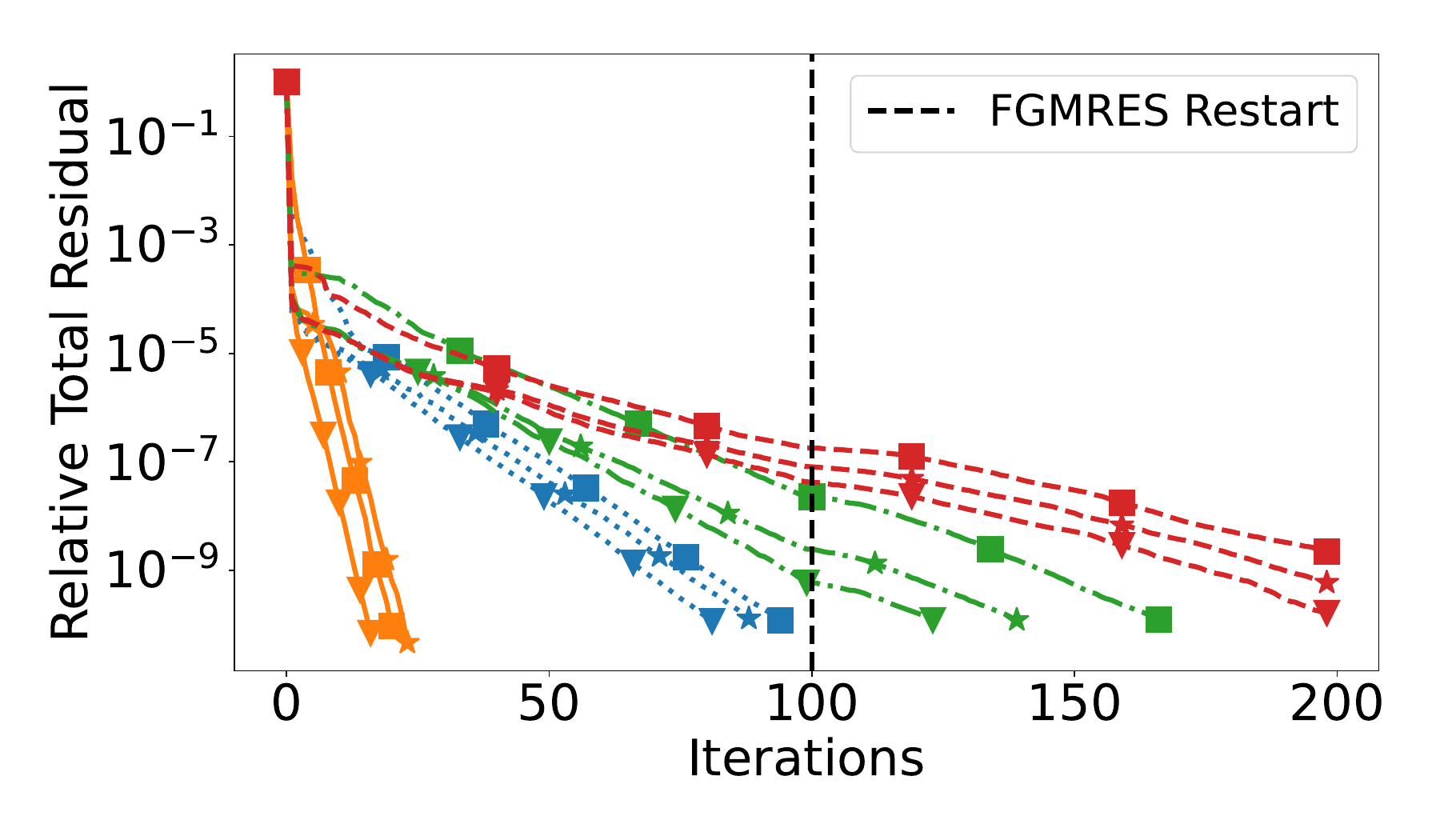}
        \includegraphics[width=0.49\textwidth]{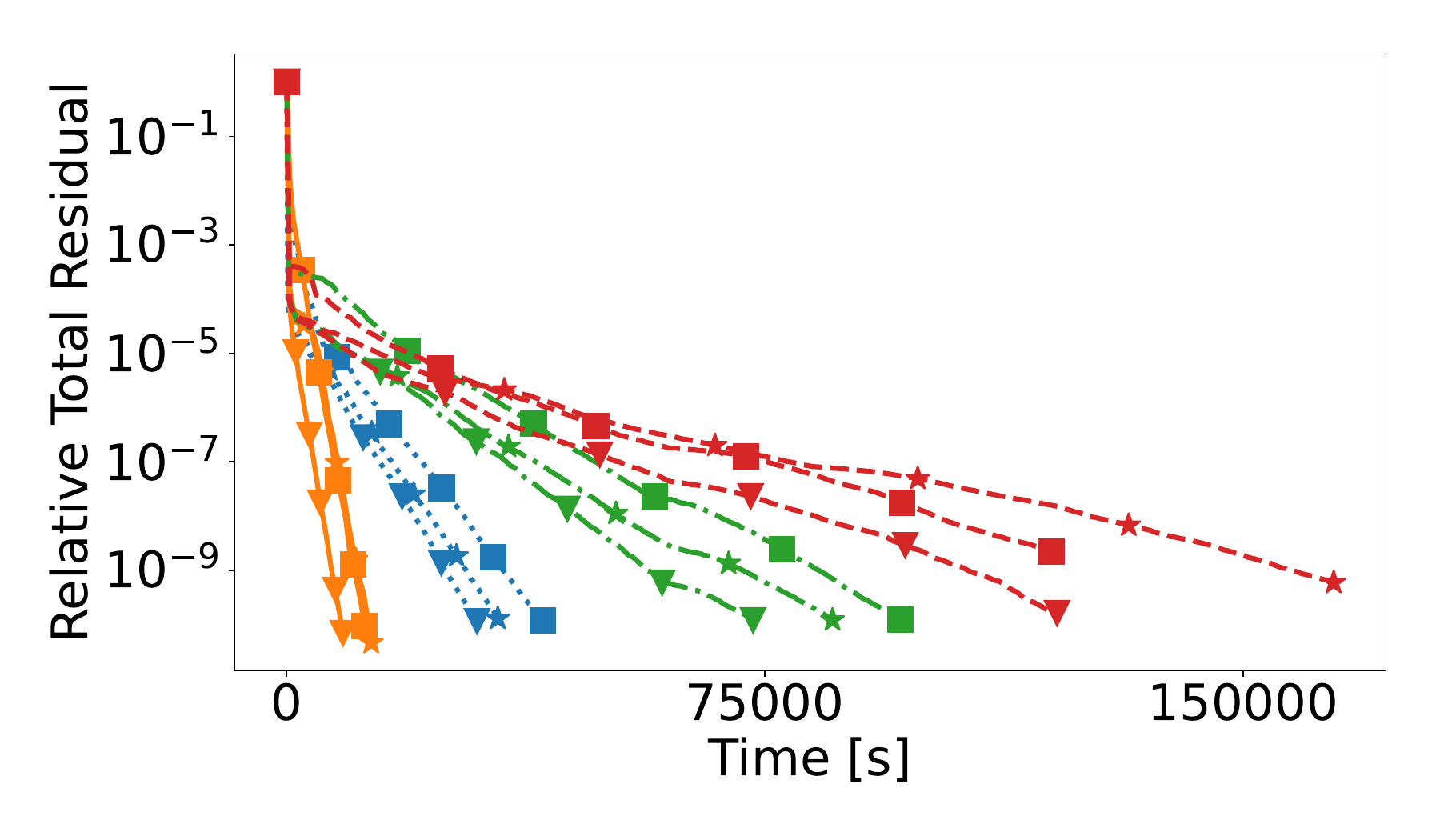}\\
        \includegraphics[width=0.49\textwidth]{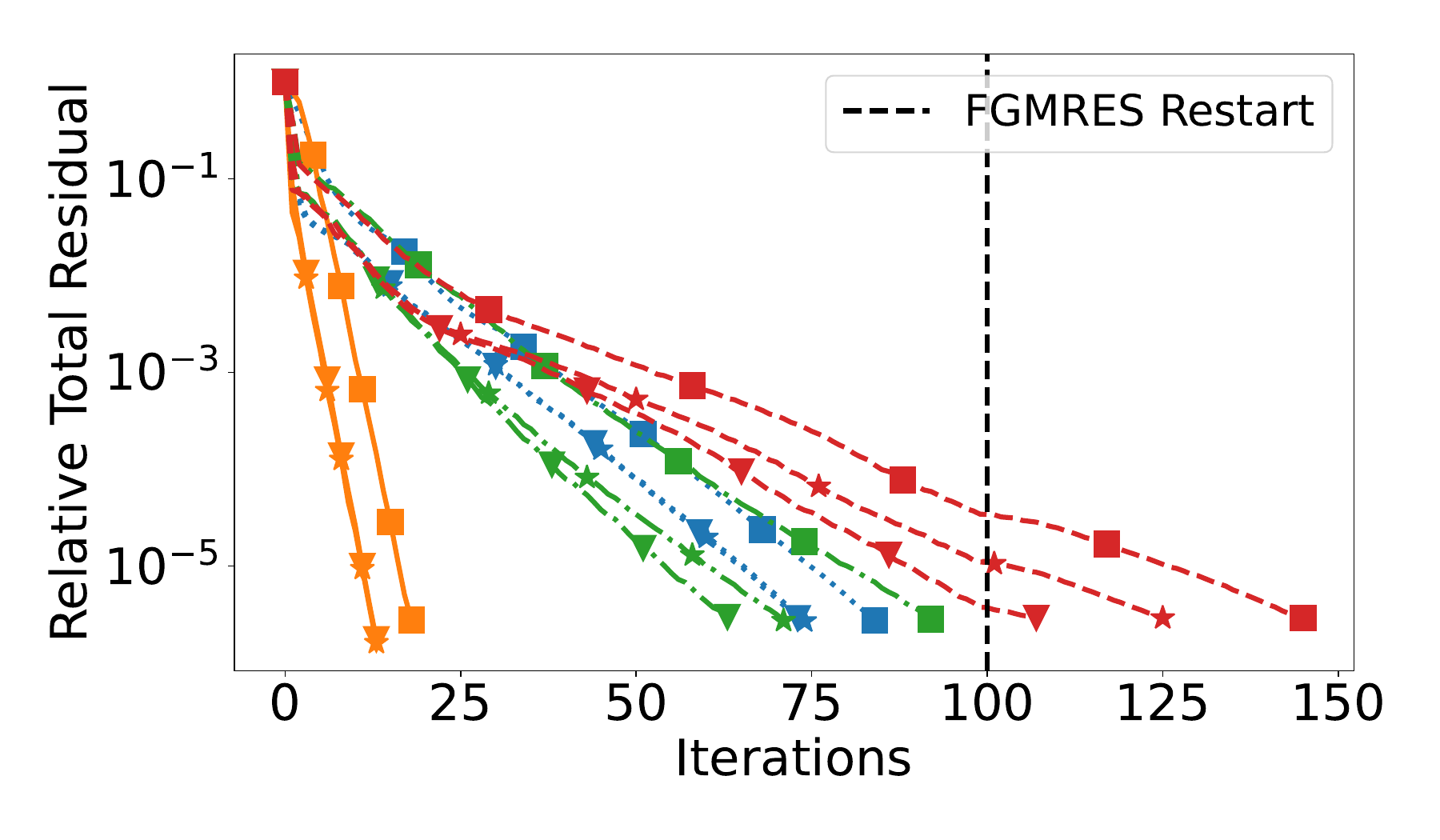}
        \includegraphics[width=0.49\textwidth]{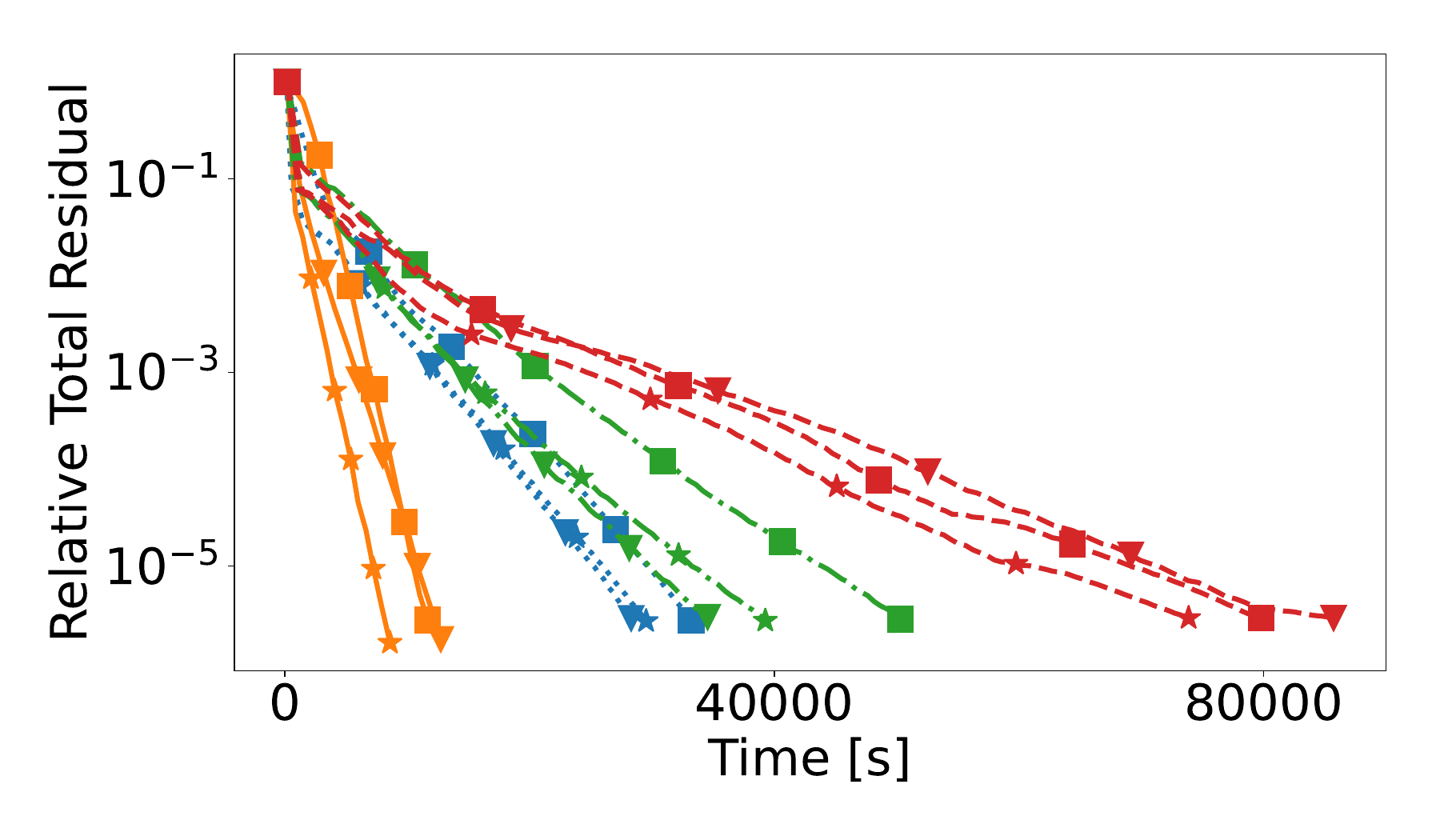}
        \\
        \includegraphics[width=0.85\textwidth]{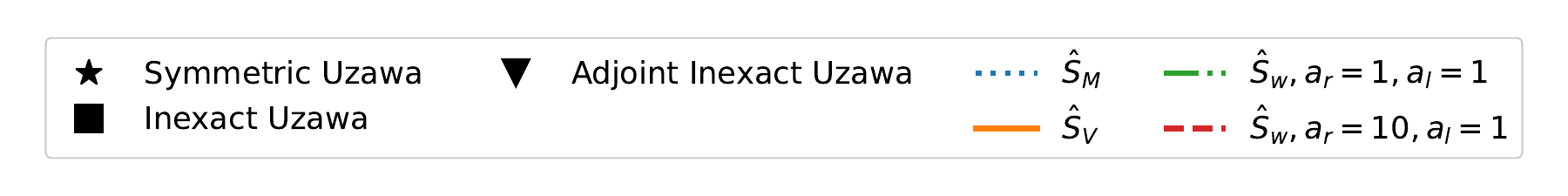}
        \caption{Saddle point relative residual over iterations (left) and time (right) with respect to the initial Stokes solve (top) or the Stokes solve after 500Myrs (bottom).}
        \label{fig:SchurCompInit}
    \end{figure}
    In this section, we compare different Schur complement approximations with the setup from Sec.~\ref{subsec:GeoModel} with $2.64 \cdot 10^8$ of degrees of freedom in $(u,p)$ and $8.45 \cdot 10^7$ degrees of freedom in $T$ based on the first Stokes solve and a Stokes solve performed after loading an exemplary temperature checkpoint at approximately 500Myrs (plate time) taken from a previous simulation using an implicit Euler time discretisation and a frozen velocity approach (compare \cite{Heister2017}) solving 
    \[
        - \nabla \cdot u^{n+1}  = \frac{\nabla \rho(x)}{\rho(x)} \cdot u_*^{n+1}
    \]
    in the saddle point system. For this test, we used a CG solver implemented in HyTeG on the coarse grid.

    Fig.~\ref{fig:SchurCompInit} depicts the relative total residual over the number of iterations and the time spent solving on a machine with 240 cores (2 Sockets with an AMD EPYC 9754 128-Core Processor @ 2.25 GHz each) for a coarse grid with 240 macro elements and a maximum refinement level $l_{\max}=6$ for all block preconditioners and Schur complement approximations mentioned in this article. In case of the inexact Uzawa and adjoint inexact Uzawa block preconditioner, $\tol_A$ has been set to $10^{-4}$, while for the symmetric Uzawa, we use $\tol_A = 10^{-2}$ to compensate for the fact that the symmetric Uzawa performs two velocity updates per application. The weighted BFBT approach is also tested with an asymmetric scaling with $a_r = 10$ and $a_l = 1$ similar to \cite{Rudi2017}.

    Fig.~\ref{fig:SchurCompInit} clearly shows that using the V-cycle BFBT approximation (see \eqref{eq:schurV}, orange, solid line) needs the lowest number of iterations and is the fastest for solving the Stokes system when reducing the relative residual by a large factor. The weighted BFBT approximation (see \eqref{eq:schurWBFBT}) with asymmetric viscosity scaling (red, dashed line) is the slowest, the mass approximation (see \eqref{eq:schurMass}, blue, dotted line) is the second fastest, in some cases comparing with the weighted BFBT approximation with symmetric viscosity scaling (green, dash-dotted line). The inexact Uzawa (see \eqref{eq:InexactUzawaStep1}, square) performs worse than the adjoint (\eqref{eq:AdjointInexactUzawaStep1}, triangle) and the symmetric Uzawa (\eqref{eq:SymmetricUzawaStep1}, star) in each of the aforementioned cases. Note that when solving to a small relative residual ($10^{-3}$ for the first time step, $10^{-1}$ for the later time step), all approximations perform comparably fast. Generally, the absolute residual is the highest for the first time step. All tests depicted in this figure, were stopped when reaching an absolute residual of $10^{-8}$. Note that the performance of the weighted BFBT heavily depends on the chosen parameters. We found that in this test case, $\omega=0.0125$ worked best. Note that the performance of the asymmetrically scaled weighted BFBT approximation seems to be problem-dependent, featuring slightly faster convergence than the symmetrically scaled weighted BFBT approximation in some of the two dimensional scenarios we investigated.
    \subsection{Scalability} \label{sec:Scaling}
    As a scaling test, we compare the first 3 time steps of a simulation with the setup from Sec.~\ref{subsec:GeoModel} with an increasing number of cores for the same problem size for the strong scaling test, and for an increasing number of degrees of freedom for the weak scaling test.
    \subsubsection{Strong Scaling}
\begin{figure}[htbp]
		\centering
        \includegraphics[width=\textwidth]{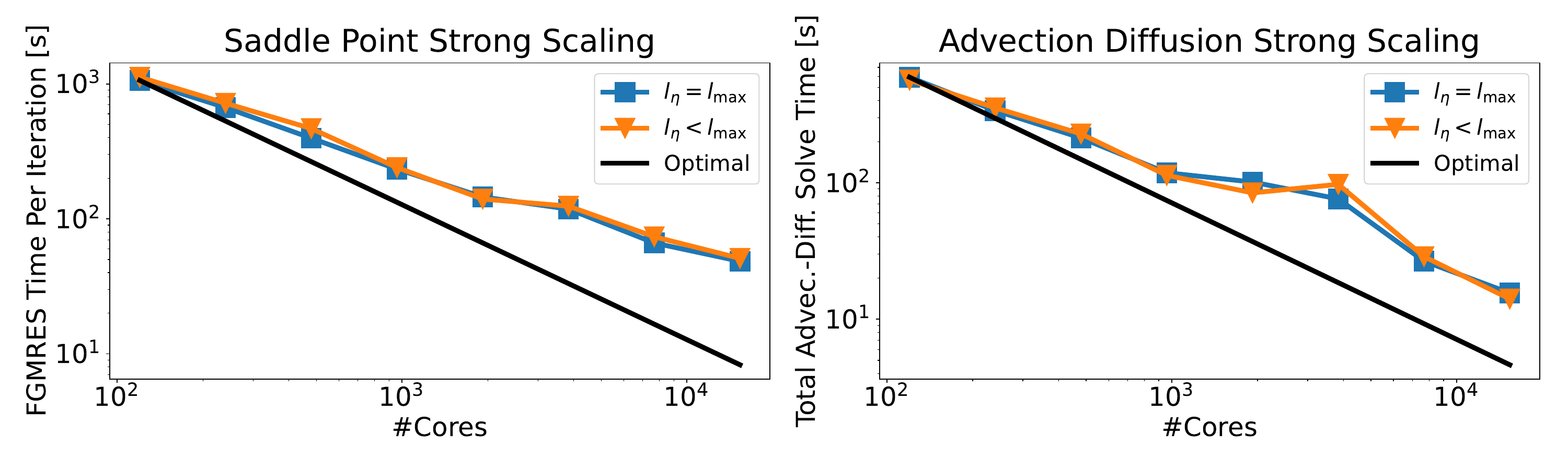}
        \\
        \begin{minipage}{0.49\textwidth}
        \begin{spacing}{0.6}
        {\footnotesize \textbf{(a)} Average time per FGMRES iteration for the saddle point problem.\textcolor{white}{This is to make the scaling of the text similar.}}
        \end{spacing}
        \end{minipage}
        \hfill
        \begin{minipage}{0.49\textwidth}
        \begin{spacing}{0.6}
        {\footnotesize \textbf{(b)} Total time for the advection-diffusion problem, consisting of FGMRES iterations for the diffusion and MMOC for the advection.}
        \end{spacing}        
        \end{minipage}        
        \caption{Strong scaling test for the saddle point problem (left) and the advection diffusion problem (right).}
        \label{fig:strongScaling}
    \end{figure}

    Given $2.64 \cdot 10^8$ of degrees of freedom in $(u,p)$ and $8.45 \cdot 10^7$ degrees of freedom in $T$, we use $120\cdot2^n$, $n=0,1,\ldots,7$ cores and measure the time needed to solve the saddle point problem per FGMRES iteration (until an absolute error of $\tol_{\rk{u,p}} =10^{-6}$ is reached) as well as the total time to solve the advection-diffusion problem, consisting of FGMRES iterations for the diffusion and MMOC for the advection. Note that in HyTeG, the maximum number of usable cores is equal to the number of macro elements in the coarse grid. The results are shown in Fig.~\ref{fig:strongScaling}. Both subproblems scale reasonably well. There is no major difference in the scaling for $l_\eta=l_{\text{max}}$ and the other option explained in Sec.~\ref{sec:ABlockTest}, marked as $l_\eta<l_{\max}$.
    \subsubsection{Weak Scaling}
\begin{figure}[htbp]
		\centering
        \begin{subfigure}{0.49\textwidth}
        \includegraphics[width=\textwidth]{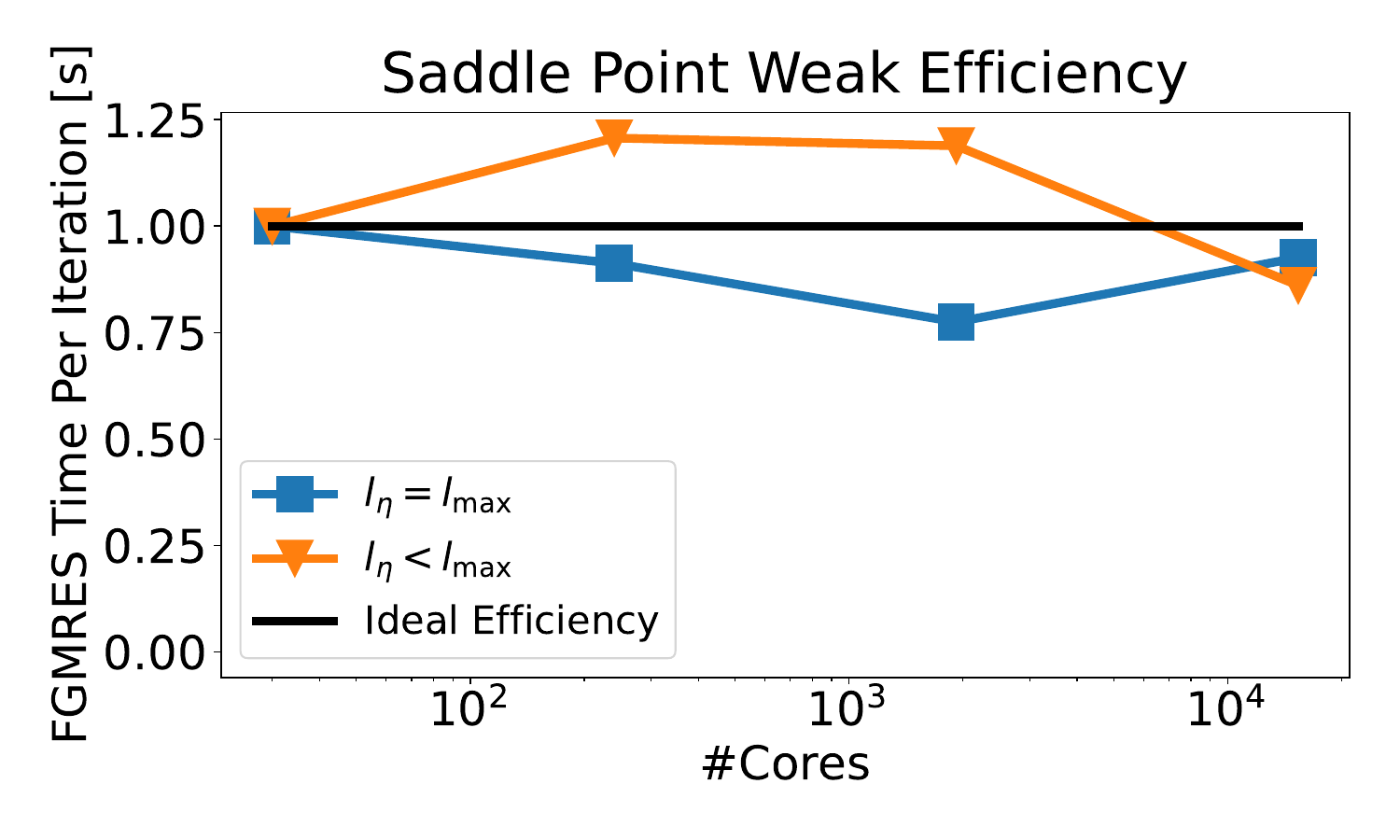}
        \caption{Weak efficiency of the average saddle point FGMRES time per iteration.}
        \label{fig:weakScalingSaddle}
        \end{subfigure}
        \begin{subfigure}{0.49\textwidth}
        \includegraphics[width=\textwidth]{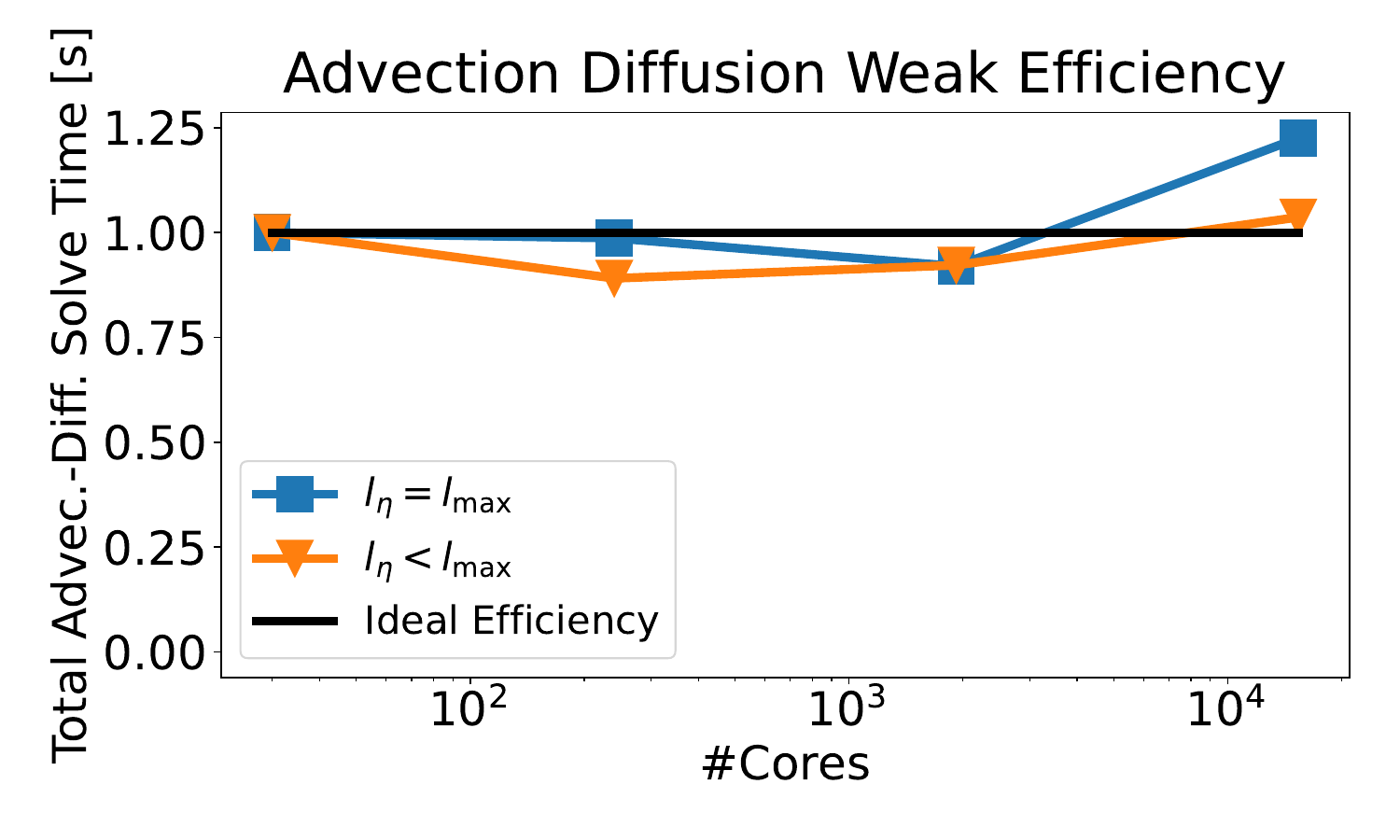}
        \caption{Weak efficiency of the total advection-diffusion FGMRES solve time.}
        \label{fig:weakScalingAdDiff}
        \end{subfigure}
        \caption{Weak scaling test for the saddle point problem (left) and the advection diffusion problem (right).}\label{fig:weakScaling}
    \end{figure}

    As a weak scaling test, we compare the first 3 time steps of a simulation with the setup from Sec.~\ref{subsec:GeoModel} with approximately $1.4 \cdot 10^5$ degrees of freedom in $(u,p)$ and $4.5 \cdot 10^4$ degrees of freedom in $T$ per core, while using $30\cdot 8^n$, $n=0,1,2,3$ cores. The Stokes system was solved up to a relative residual tolerance of $10^{-6}$, setting $\tol_{\rk{u,p}}$ to $10^{-14}$ which is not reached in this test, and a fixed time step size of $0.125$ Myrs. Note that for $l_{\eta} < l_{\max}$, the evaluation of the viscosity at the quadrature points is performed exclusively as an evaluation of $\eta^{n+1}_{\text{P1}}$ for $15\,360$ cores. Fig.~\ref{fig:weakScaling} depicts the results. As with the strong scaling test, both subproblems scale reasonably well, sometimes better than expected: In the saddle point problem, the code scales better than ideal if $l_\eta<l_{\text{max}}$ and going from 30 to 240 cores, and in the advection-diffusion problem, the code scales much better than expected if $l_\eta=l_{\text{max}}$ going from $1\,920$ cores to $15\,360$ cores. In both cases, although not definitive, we suspect that the higher resolution in space yields a smoother viscosity in the shear heating and the $A$ block, resulting in a simpler problem to solve. Note that as previously described, the initial state has a fixed random component in each core, therefore small changes in the run-time may occur.
  \section{Conclusion} \label{sec:Conclusion}
  In this study, we investigated numerical schemes to solve the TALA, a common approach to model the convective flow in the Earth's mantle. While approximating the spatial structure via blending, the proposed approach uses a splitting method between the time dependent equation for the temperature and the time independent Stokes system for the velocity and the pressure. Our approach combines a BDF2 scheme with variable time step resulting in a second order convergence in time. While a frozen velocity approach for tackling the compressible case is quite popular, we use the physically more accurate implicit version of the compressible mass balance. Although this results in a non-symmetric off-diagonal block in the saddle point Stokes formulation, our iterative solver is robust with respect to this term. The overall code scales well enough such that high-resolution simulations are possible. Solving the saddle point problem with $2.11\cdot10^9$ unknowns while having $6.7\cdot10^8$ unknowns in the temperature needs 69\% of the overall computation time. Therefore, different approaches to solve the saddle point problem have been investigated. They showed that when only a small error reduction is needed per time step, all solvers perform well, while a larger reduction is fastest reached with the V-cycle BFBT approximation. High-resolution simulations with $2.7\cdot10^9$ degrees of freedom in space and 6166 steps in time showed the effect of plate movement on the creation of slabs, plumes, and convection cells.
\section*{Acknowledgments} \label{sec:Acknowledgments}
We thank Berta Vilacis for her advice on selecting suitable viscosity models and her work on extracting and slightly modifying the plate reconstructions from \cite{Mueller2022}\footnote{model data freely available at \href{https://www.earthbyte.org/webdav/ftp/Data_Collections/Muller_etal_2022_SE/}{https://www.earthbyte.org/webdav/ftp/Data\_Collections/Muller\_etal\_2022\_SE/}} using GPlates \cite{Mueller2018}\footnote{software freely available at \href{https://www.gplates.org/}{https://www.gplates.org/} } as well as making them available in HyTeG. We further want to thank Andreas Wagner and Marcus Mohr for helpful discussions about HyTeG. 

\section*{Declarations}
\subsection*{Funding}
Funding for this work was provided by the German Research Foundation (DFG) as part of the research training group 440512084 UPLIFT\footnote{\href{https://gepris.dfg.de/gepris/projekt/440512084}{https://gepris.dfg.de/gepris/projekt/440512084}} (GRK 2698) and the CoMPS: Multi-physics simulations for Geodynamics on heterogeneous Exascale Systems project funded by the German Federal Ministry of Education and Research (BMBF) via its SCALEXA (New Methods and Technologies for Exascale Computing) initiative.

The numerical calculations presented in this article were carried out on the HAWK Supercomputer\footnote{\href{https://www.hlrs.de/de/loesungen/systeme/hpe-apollo-hawk}{https://www.hlrs.de/de/loesungen/systeme/hpe-apollo-hawk}} at the H\"ochstleistungrechenzentrum Stuttgart (HLRS) as part of the project CoMPSatHLRS (ACID 44279) and on local machines at the Technische Universit\"at M\"unchen (TUM) Department of Mathematics. During development, test calculations were also carried out on the TETHYS-3G\footnote{\href{https://www.geophysik.uni-muenchen.de/en/research/geocomputing/facilities/tethys-3g}{https://www.geophysik.uni-muenchen.de/en/research/geocomputing/facilities/tethys-3g}} Cluster \cite{Oeser2006}. TETHYS-3G is being co-funded by the Free State of Bavaria and the Deutsche Forschungsgemeinschaft (DFG, German Research Foundation) in the framework of the programme for Forschungsgroßgeräte nach Art. 91b GG – 495931446.
\subsection*{Conflict of Interest}
The authors have no competing interests to declare that are relevant to the content of this article.

\begin{appendices}
\section{Source Code}
The full source code of the implementations presented in this article as well as usage instructions are available at \url{https://www.doi.org/10.5281/zenodo.15497635}.

\section{Solver Parameters} \label{solverParameters}
\vspace{-1.0cm}
 \begin{table}[h] 
     \caption{\label{parametertable} Default solver parameters chosen for the simulations described in Sec.~\ref{sec:GeodynamicApplication}, unless stated otherwise.} 
    \begin{center} 
        \text{ }
        \begin{tabular}{|p{0.12\linewidth}|p{0.06\linewidth}|p{0.75\linewidth}|}
        \hline
        Symbol & Value  & Description \\
        \hline
        $\sigma$ & $1$ & Scaling of $\hat{A}^{-1}$ in block preconditioners \\
        $\omega$ & $0.3$ & Scaling of $\hat{S}^{-1}$ in block preconditioners \\
        $\tol_{A}$ & $10^{-2}$ & Relative tolerance of the outer CG $\hat{A}$ solver \\        
        $m_{A}$ & $3$ & Number of Chebyshev pre- and postsmoothing steps applied on each refinement level while performing a V-cycle w.r.t. $\hat{A}$ \\     
        $\deg_{A}$ & $2$ & Polynomial degree of the Chebyshev smoother w.r.t. $\hat{A}$  \\  
        $\tol_{\VBFBT}$ & $10^{-1}$ & Relative tolerance of the outer CG $\rk{B\hat{A}_C^{-1} B^T}$ solver used while applying the V-cycle BFBT Schur complement approximation  \\  
        $\tol_{\invMass}$ & $10^{-10}$ & Relative tolerance of the outer CG $\hat{S}_M$ solver used while applying the inverse viscosity scaled mass Schur complement approximation and while preconditioning $\rk{B\hat{A}_C^{-1} B^T}$ with $\hat{S}_M$ as part of the V-cycle BFBT Schur complement approximation  \\
        $\tol_{\wBFBT}$ & $10^{-10}$ & Relative or absolute tolerance (whichever applies first) of the outer CG $K_{1/\sqrt{\eta}}$ solver used while applying the w-BFBT Schur complement approximation  \\        
        $m_V$ & $1$ & Number of Chebyshev pre- and postsmoothing steps applied on each refinement level while performing a V-cycle w.r.t. $\hat{A}_C$  \\    
        $\deg_V$ & $1$ & Polynomial degree of the Chebyshev smoother w.r.t. $\hat{A}_C$  \\         
        $\tol_{\coarse}$ & $10^{-2}$ & Relative tolerance of the CG coarse grid solver w.r.t. $\hat{A}$ and $\hat{A}_C$  \\  
        $\tol_{\rk{u,p}}$ & $10^{-5}$ & Absolute tolerance of the saddle point FGMRES solver  \\      
        $\tol_{T}$ & $10^{-10}$ & Absolute tolerance of the temperature point FGMRES solver  \\  
        $\tol_{\VectorMass}$ & $10^{-4}$ & Relative tolerance of the outer CG ${M}_{\sqrt{\eta_r}}^{-1}$ and ${M}_{\sqrt{\eta_l}}$ solver \\          
        $a_r, a_l$ & $1$ & Potentially asymmetric w-BFBT scaling of ${M}_{\sqrt{\eta_r}}^{-1}$,${M}_{\sqrt{\eta_l}}$,$K_{1/\sqrt{\eta_l}}^{-1}$ and $K_{1/\sqrt{\eta_r}}^{-1}$  \\

        \hline
    \end{tabular}
    \end{center} 
\end{table}
\newpage    
	\section{Nondimensionalisation\label{Nondim}}
\subsection{Nondimensionalisation of the TALA}
		We choose a nondimensionalisation of the TALA similar to \cite[Chap. 6.10]{Schubert2001} and \cite{Leng2008}. Tab.~\ref{dimtable} shows the reference constants used for the nondimensionalisation of the forward mantle convection simulation presented in Sec.~\ref{subsec:GeoModel}. In Sec.~\ref{subsec:ChoiceOfReferenceConstants}, we explain how the specific reference values in Tab.~\ref{dimtable} were chosen. 

		Let $i \in \left \{1, \ldots, \dim \right \}$. First, we choose reference constants and introduce nondimensional variables for the space, velocity, temperature, density, gravity, viscosity, thermal expansivity, Grüneisen parameter, specific heat capacity and thermal conductivity by
		\begin{align*}
			x_i 	&=: d \tilde{x}_i, & \rho &=: \rho_0 \tilde{\rho}, & \Gamma &=: \Gamma_0 \tilde{\Gamma},	\\
			u 	&=: u_0 \tilde{u},	& g &=: g_0 \tilde{g},	& C^p &=: C^p_0 \tilde{C}^p,		\\
			T 	&=: \Delta T \tilde{T}, & \eta	&=: \eta_0 \tilde{\eta}, & k &=: k_0 \tilde{k},			\\
			T_d &=: \Delta T \tilde{T}_d,	&  \alpha 	&=: \alpha_0 \tilde{\alpha}. 							& 			&
		\end{align*}
        The dimensionless Grüneisen parameter is defined as $\Gamma := \frac{\alpha K^S}{\rho C^p}$ with $K^S$ denoting the isentropic bulk modulus. From this set of chosen reference constants, we derive nondimensionalisation constants for the time, internal heating, pressure and thermal diffusivity via
		\begin{align*}
			t 	&=: {\color{black} \frac{d}{u_0}} \tilde{t} =: {\color{black} t_0} \tilde{t},	& p 	&=: {\color{black} \frac{\eta_0 u_0}{d}} \tilde{p} =: {\color{black} p_0} \tilde{p}, 	& {\color{black} \kappa^d_0} 	:= {\color{black} \frac{k_0}{\rho_0 C^p_0}},	\\
			H 	&=: {\color{black} \frac{C^p_0 \Delta T u_0}{d}} \tilde{H} =: {\color{black} H_0} \tilde{H}, & p_d &=: {\color{black} \frac{\eta_0 u_0}{d}} \tilde{p}_d =: {\color{black} p_0} \tilde{p}_d. 	  	& 
		\end{align*}
	      Now we can introduce the dimensionless Péclet number, Rayleigh number and Dissipation number as
		\begin{align*}
			\Peclet	&:= \frac{u_0 d}{\kappa^d_0},	& \Rayleigh	&:= \frac{\rho_0 \alpha_0 g_0 \Delta T d^3}{\kappa^d_0 \eta_0}, 	& \Dissipation &:= \frac{\alpha_0 g_0 d}{C^p_0}.
		\end{align*}
        For the sake of completeness, we also introduce nondimensionalisation constants for the isothermal compressibility $\kappa^T$, isothermal bulk modulus $K^T$ and latent heat generated by phase transitions $Q^L$ (compare \cite{Gassmoeller2020}) via

		\begin{align*}
			\kappa^T &=: {\color{black} \frac{\alpha_0}{\Gamma_0 \rho_0 C^p_0}} \tilde{\kappa}^T =: {\color{black} \kappa^T_0} \tilde{\kappa}^T,
            & K_0^T &:= \frac{1}{\kappa_0^T}, 	
            & Q^L 	&=: {\color{black} \frac{\rho_0 u_0 \Delta T C^p_0}{d}} \tilde{Q}^L =: {\color{black} Q^L_0} \tilde{Q}^L,
		\end{align*}
        and define the dimensionless mantle compressibility $\gamma$ and driving term $\xi$ (compare \cite{Ricard2007} by
        \begin{align*}
			\gamma	:= \frac{\Dissipation}{\Gamma_0}, \quad \xi := \alpha_0 \Delta T.	
		\end{align*}

    		\subsection{Choice of Nondimensionalisation Constants} \label{subsec:ChoiceOfReferenceConstants}
            \vspace{-0.8cm}
        \begin{table}[h] 
            \text{ }
            \begin{center} 
            \text{ }
            \caption{\label{dimtable} Reference constants for nondimensionalisation and specific model values.} 
                \begin{tabular}{|l|c|l|l|}
					\multicolumn{4}{c}{Independent reference constants} \\
					\hline
					Name & Symbol & Value & Unit \\
					\hline
					Mantle thickness & $d$ & $2.891 \cdot 10^{6}$ & $\um$ \\
					Temperature difference & $\Delta T$ & $3.900 \cdot 10^{3}$ & $\uK$ \\
					Viscosity & $\eta_0$ & $1.000 \cdot 10^{22}$ & $\uPa \us$ \\
					Density & $\rho_0$ & $4.686 \cdot 10^{3}$ & $\ukg \um^{-3}$ \\
					Specific heat capacity & $C^p_0$ & $1.250 \cdot 10^{3}$ & $\uJ \uK^{-1} \ukg^{-1}$ \\
					Thermal expansivity & $\alpha_0$ & $2.000 \cdot 10^{-5}$ & $\uK^{-1}$ \\
					Gravity & $g_0$ & $9.810$ & $\um \us^{-2}$ \\
					Velocity & $u_0$ & $5.000 \cdot 10^{-9}$ & $\um \us^{-1}$ \\
					Thermal conductivity & $k_0$ & $3.000$ & $\uW \um^{-1} \uK^{-1}$ \\
					Grüneisen parameter & $\Gamma_0$ & $1.200$ &  \\
					\hline
					\multicolumn{4}{c}{ } \\
					\multicolumn{4}{c}{Dependent reference constants} \\
					\hline
					Name & Symbol & Value & Unit \\
					\hline
					Thermal diffusivity & $\kappa^d_0$ & $5.122 \cdot 10^{-7}$ & $\um^2 \us^{-1}$ \\
					Time & $t_0$ & $5.782 \cdot 10^{14}$ & $\us$ \\
					Pressure & $p_0$ & $1.730 \cdot 10^{7}$ & $\uPa$ \\
					Internal heating & $H_0$ & $8.431 \cdot 10^{-9}$ & $\uW \ukg^{-1}$ \\
					Latent heat generated by phase Transitions & $Q^L_0$ & $3.951 \cdot 10^{-5}$ & $\uW \um^{-3}$ \\
					Isothermal bulk modulus & $K^T_0$ & $3.514 \cdot 10^{11}$ & $\uPa$ \\
					Isothermal compressibility & $\kappa^T_0$ & $2.845 \cdot 10^{-12}$ & $\uPa^{-1}$ \\
					Rayleigh number & $\Rayleigh$ & $1.692 \cdot 10^{7}$ &  \\
					Dissipation number & $\Dissipation$ & $4.538 \cdot 10^{-1}$ &  \\
					Péclet number & $\Peclet$ & $2.822 \cdot 10^{4}$ &  \\
					Mantle compressibility & $\gamma$ & $3.781 \cdot 10^{-1}$ &  \\
					Driving term & $\xi$ & $7.800 \cdot 10^{-2}$ &  \\
					\hline
					\multicolumn{4}{c}{ } \\
					\multicolumn{4}{c}{Specific Values} \\
					\hline
					Name & Symbol & Value & Unit \\
					\hline
					Radius surface & $r_{\Surface}$ & $6.371 \cdot 10^{6}$ & $\um$ \\
					Radius CMB & $r_{\CMB}$ & $3.480 \cdot 10^{6}$ & $\um$ \\
					Surface temperature & $T_{\Surface}$ & $3.000 \cdot 10^{2}$ & $\uK$ \\
					CMB temperature & $T_{\CMB}$ & $4.200 \cdot 10^{3}$ & $\uK$ \\
					Density at the top of the mantle & $\rho_{\Top}$ & $3.381 \cdot 10^{3}$ & $\ukg \um^{-3}$ \\
					Adiabatic temperature at the top of the mantle & $T_{\text{adiabatic}}$ & $1.6 \cdot 10^{3}$ & $\uK$ \\
					\hline
			\end{tabular}
            \end{center} 
        \end{table}
		The surface radius $r_{\Surface}$, CMB radius $r_{\CMB}$ and resulting mantle thickness $d$ are taken directly from the preliminary reference Earth model (PREM) \cite{Dziewonski1981}. The density $\rho_{\Top}$ matches the density value at the top of the upper mantle at radius 6346.6$\ukm$ and the density reference constant $\rho_0$ matches the average density in the mantle (3480$\ukm$ - 6346.6$\ukm$) given by the PREM. The approximate gravitational acceleration $9.81 \frac{\um}{\us^2}$ at the surface of the Earth is chosen as $g_0$.

        As mentioned in Sec.~\ref{sec:Introduction}, different estimates for the temperature at the CMB are available ranging from $2700\uK\pm250\uK$ \cite{Jeanloz1986} to $4500\uK\pm 1000\uK$ \cite{daSilva2000, Jeanloz1990, Alfe1999} (also compare the summaries given in \cite[Chapter 7.7]{Fowler2004} and \cite{Bukowinski1999,Knittle1988}). We've chosen a CMB temperature of $4200\uK$ for our simulations which approximately matches the pyrolite estimations of 4250$\uK\pm250\uK$ in \cite{daSilva2000}.

		The dimensionless Grüneisen parameter $\Gamma$ is estimated to lie between 1.1 and 1.3 for the Earth's lower mantle \cite{Stacey2004} and between 1.4 and 1.5 for the outer core \cite{Alfe2002, Stacey2004}. We've chosen $\Gamma_0 = 1.2$ as our reference value (compare \cite{Colli2020} using $\Gamma_0 = 1.1$).

		For the lower mantle the thermal conductivity $k$ is assumed to be in the range of $2.5 \frac{\uW}{\um \uK}$ to $3.5 \frac{\uW}{\um \uK}$ \cite{Tang2014} (compare the values in \cite{Ricard2007, Busse1994, Colli2020}). Whilst it is feasible that the thermal conductivity in upper mantle could be upwards of $5.5 \frac{\uW}{\um \uK}$, see \cite{Grose2019}, we still choose $k_0 = 3.0 \frac{\uW}{\um \uK}$ for our reference constant as the lower mantle makes up about two thirds of the total mantle volume. 
		Larger thermal conductivity values in the range of $4.0 \frac{\uW}{\um \uK}$ to $4.7 \frac{\uW}{\um \uK}$ are also commonly used (compare \cite{Liu2019, Gassmoeller2020, Dannberg2016}).

        Estimations for the thermal expansivity $\alpha$ of the mantle lie between $1 \cdot 10^{-5} \frac{1}{\uK}$ and $ 3.0 \cdot 10^{-5} \frac{1}{\uK}$ \cite{Stacey2004} (compare the values used in \cite{Ricard2007, Blankenbach1989, Busse1994, Gassmoeller2020, Liu2019, Dannberg2016, Colli2020}. We've chosen $\alpha_0 = 2.0 \cdot 10^{-5} \frac{1}{\uK}$ as suggested in \cite{Ricard2007, Gassmoeller2020, Dannberg2016}.

		Commonly used values for the specific heat capacity $C^p_0$ of the mantle are in the range of 1000$\frac{\uJ}{\uK \ukg}$ to 1250$\frac{\uJ}{\uK \ukg}$ (compare \cite{Ricard2007, Gassmoeller2020, Blankenbach1989, Busse1994, Liu2019, Colli2020}). Similar to \cite{Gassmoeller2020} and \cite{Blankenbach1989} we set $C^p_0 = 1250$$\frac{\uJ}{\uK \ukg}$.

\begin{figure}[htbp]
    		\centering
            \includegraphics[width=\textwidth]{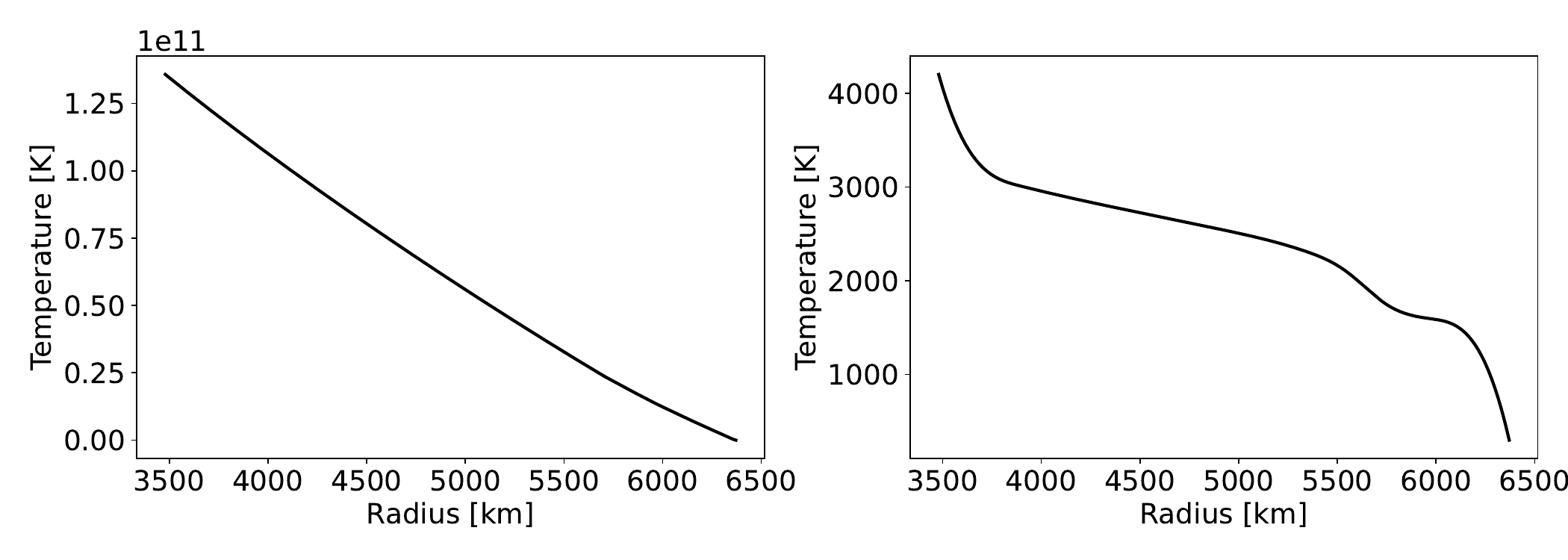}
            \caption{ Depth dependent pressure from the PREM \cite{Dziewonski1981} (left) and a temperature profile between 300$\uK$ and 4200$\uK$ qualitatively modeled after \cite[Figure 7.16]{Fowler2004} (right).}
            \label{fig:depthProfiles}
    	\end{figure}

        Our values for the thermal expansivity, Grüneisen parameter and specific heat capacity have been verified using the mineralogical thermodynamics framework MMA-EoS\footnote{Software freely available at \url{https://www.chust.org/repos/eos} under a GNU General Public License Version 3.} \cite{Chust2017} together with the thermodynamics dataset from \cite{Stixrude2011} and the assumption of a NCFMAS pyrolite bulk chemical composition (compare \cite[Sec.~5]{Chust2017}) for the Earth’s mantle. However, the validity of this choice of chemical composition especially for the lower mantle is an open research topic (compare \cite[Sec.~5.1.1]{Chust2017}, \cite{Murakami2012}). The framework allows the user to calculate (among other properties) $\rho$, $C^p$, $\alpha$, $\kappa^T$ and $\Gamma$ as a function of pressure and temperature by means of finding an optimal phase assemblage through minimizing the Gibbs-Energy (compare \cite{Chust2017}). The density $\rho$, isothermal compressibility $\kappa^T$, isobaric thermal expansivity $\alpha^p$, the Grüneisen parameter $\Gamma$ and the isobaric specific heat capacity $C^p$ are averaged along a realistic depth dependent pressure and temperature reference profile of the Earth's mantle, depicted in Fig.~\ref{fig:depthProfiles}. The average values are shown in Tab.~\ref{avgtable}, closely matching our chosen reference values (compare \cite[Sec.~5.3]{Weismueller2015}). 
        \begin{table}[h]
            \centering
            \caption{Thermodynamic properties and their mean along the depth dependent profiles depicted in Fig.~\ref{fig:depthProfiles}. \label{avgtable}}
            \begin{tabular}{|l|l|}
                \hline
                Parameter & Mean \\
                \hline
                $\rho$ & $4.645 \cdot 10^{3}$ \\
                $\kappa^T$ & $3.210 \cdot 10^{-12}$ \\
                $\alpha^p$ & $1.908 \cdot 10^{-5}$ \\
                $C^p$ & $1.233 \cdot 10^{3}$ \\
                $\Gamma$ & $1.233$ \\
                \hline
            \end{tabular}
            
	   \end{table}
       \noindent
        Estimates of the Earth’s mantle viscosity range from $10^{19}$$\uPa \us$ in the asthenosphere to $10^{23}$$\uPa \us$ in the lower mantle \cite[Chapter 5.7]{Fowler2004} hence a wide range of viscosity reference values $\eta_0$ can be found in literature (compare \cite{Ricard2007, Gassmoeller2020, Blankenbach1989, Busse1994, Liu2019, Dannberg2016, Colli2020}).
        As the remaining driving factors for the nondimensional Péclet number and Rayleigh number we've chosen $\eta_0 = 1.0 \cdot 10^{22} \uPa \us$ and $u_0 = 5.0 \cdot 10^{-9} \frac{\um}{\us}$ to get $\Peclet \approx 10^4$ (compare \cite{Colli2020}, \cite[Sec.~2.1.2]{Weismueller2015}) and $\Rayleigh \approx 10^7$ (compare \cite{Ricard2007}).    

        \subsection{Nondimensionalisation of the Governing Equations} \label{subsec:NondimEq}
        Starting with a rescaling of the spatial and temporal variables in \ref{subsec:spatialAndTemporalDerivatives}, we provide the nondimensional form of momentum, mass and energy balance equations in \ref{subsec:momEqu}, \ref{subsec:massEqu} and \ref{subsec:energyEqu}, respectively.
		\subsubsection{Spatial and Temporal Derivatives}\label{subsec:spatialAndTemporalDerivatives}
		\begin{align}
			\frac{\partial}{\partial x_i} &= \frac{\partial \tilde{x}_i}{\partial x_i} \frac{\partial}{\partial \tilde{x}_i} = \frac{1}{d} \frac{\partial}{\partial \tilde{x}_i} \quad \Rightarrow \quad \nabla = \frac{1}{d} \tilde{\nabla} := \frac{1}{d} \begin{pmatrix} \frac{\partial}{\partial \tilde{x}_1}\\ \vdots \\ \frac{\partial}{\partial \tilde{x}_{\dim}} \end{pmatrix} 
            
            \nonumber \\
            
			\frac{\partial}{\partial t} &= \frac{\partial \tilde{t}}{\partial t} \frac{\partial}{\partial \tilde{t}} = \frac{u_0}{d} \frac{\partial}{\partial \tilde{t}} \label{eq:tildeDerivatives}
		\end{align}
		\subsubsection{Momentum Conservation Equation}\label{subsec:momEqu}
		We start with the dimensional momentum conservation equation
		\begin{align}
			- \nabla \cdot \tau + \nabla p_d + \rho \alpha T_d g = 0.
		\end{align}
        Using \eqref{eq:tildeDerivatives}, we get
		\begin{align*}
			- \nabla \cdot \tau &= - \nabla \cdot \eta \left (  \nabla u + \left ( \nabla u \right )^T - \frac{2}{\dim} \left ( \nabla \cdot u \right ) I  \right )
			\\
			&= -\frac{u_0 \eta_0}{d^2} \tilde{\nabla} \cdot \tilde{\eta} \left ( \tilde{\nabla} \tilde{u} + \left ( \tilde{\nabla} \tilde{u} \right )^T - \frac{2}{\dim} \left ( \tilde{\nabla} \cdot \tilde{u} \right ) I  \right )
			\\
			&=: -\frac{u_0 \eta_0}{d^2} \tilde{\nabla} \cdot \tilde{\tau},
		\end{align*}
        \vspace{-0.5cm}
		\begin{align*}
			\nabla p_d &= \frac{u_0 \eta_0}{d^2} \tilde{\nabla} \tilde{p}_d, & \rho \alpha T^d g &= \frac{u_0 \eta_0}{d^2} \frac{\Rayleigh}{\Peclet} \tilde{\rho} \tilde{\alpha} \tilde{T}^d \tilde{g}.
		\end{align*}
		After dividing by $\frac{u_0 \eta_0}{d^2}$ and dropping the tilde in notation, the nondimensional momentum conservation equation reads
		\begin{align}
			- \nabla \cdot \tau + \nabla p_d + \frac{\Rayleigh}{\Peclet} \rho \alpha T^d g = 0.
		\end{align}
		\subsubsection{Mass Conservation Equation} \label{subsec:massEqu}
		We start with the dimensional mass conservation equation
		\begin{align}
			 - \nabla \cdot u - \frac{\nabla \rho}{\rho} \cdot u = 0.
		\end{align}
        Using \eqref{eq:tildeDerivatives}, we get
		\begin{align*}
			 - \nabla \cdot u = - \frac{u_0}{d} \tilde{\nabla} \cdot \tilde{u},
             \quad \quad \quad 
        - \frac{\nabla \rho}{\rho} \cdot u = - \frac{u_0}{d} \frac{\tilde{\nabla} \tilde{\rho}}{\tilde{\rho}} \cdot \tilde{u}.
		\end{align*}
		After dividing by $\frac{u_0}{d}$ and dropping the tilde in notation the nondimensional mass conservation equation reads
		\begin{align}
			 - \nabla \cdot u - \frac{\nabla \rho}{\rho} \cdot u = 0.
		\end{align}
		\subsubsection{Energy Conservation Equation}\label{subsec:energyEqu}
		We start with the dimensional energy conservation equation
		\begin{align}
			\rho C^p \left ( \frac{\partial T}{\partial t} + u \cdot \nabla T \right ) - \nabla \cdot \left ( k \nabla T \right ) - \alpha T \rho \left ( u \cdot g \right )  - \rho H - \tau : \varepsilon  = 0.
		\end{align}
        Using \eqref{eq:tildeDerivatives}, we get
		\begin{align*}
			\rho C^p \left ( \frac{\partial T}{\partial t} + u \cdot \nabla T \right ) = \frac{\rho_0 C^p_0 \Delta T u_0}{d} \tilde{\rho} \tilde{C}^p \left ( \frac{\partial \tilde{T}}{\partial \tilde{t}} + \tilde{u} \cdot \tilde{\nabla} \tilde{T} \right ),
		\end{align*}
		\begin{align*}
			- \nabla \cdot \left ( k \nabla T \right ) = -\frac{k_0 \Delta T}{d^2} \tilde{\nabla} \cdot \left ( \tilde{k} \tilde{\nabla} \tilde{T} \right ) = - \frac{\rho_0 C^p_0 \Delta T u_0}{d}  \frac{1}{\Peclet} \tilde{\nabla} \cdot \left ( \tilde{k} \tilde{\nabla} \tilde{T} \right ),
		\end{align*}
		\begin{align*}
			- \tau : \varepsilon = - \frac{\eta_0 u_0^2}{d^2} \tilde{\tau} : \tilde{\varepsilon} = - \frac{\rho_0 C^p_0 \Delta T u_0}{d} \frac{\Peclet \Dissipation}{\Rayleigh} \tilde{\tau} : \tilde{\varepsilon},
		\end{align*}
		\begin{align*}
			- \alpha T \rho \left ( u \cdot g \right ) = - \alpha_0 \Delta T \rho_0 u_0 \tilde{\alpha} \tilde{T} \tilde{\rho} \left ( \tilde{u} \cdot \tilde{g} \right ) = - \frac{\rho_0 C^p_0 \Delta T u_0}{d} \Dissipation \tilde{\alpha} \tilde{T} \tilde{\rho} \left ( \tilde{u} \cdot \tilde{g} \right ),
		\end{align*}
		\begin{align*}
			- \rho H = -  \rho_0 H_0 \tilde{\rho} \tilde{H} = - \frac{\rho_0 C^p_0 \Delta T u_0}{d}  \tilde{\rho} \tilde{H}.
		\end{align*}
		After dividing by $\frac{\rho_0 C^p_0 \Delta T u_0}{d}$ and dropping the tilde in notation, the nondimensional energy conservation equation reads
		\begin{align}
			\rho C^p \left ( \frac{\partial T}{\partial t} + u \cdot \nabla T \right ) - \frac{1}{\Peclet} \nabla \cdot \left ( k \nabla T \right ) - \Dissipation \alpha T \rho \left ( u \cdot g \right ) - \rho H - \frac{\Peclet \Dissipation}{\Rayleigh} \tau : \varepsilon  = 0.	
		\end{align}
	\end{appendices}

	
	\bibliography{sn-bibliography}

\end{document}